\documentstyle[fullpage,12pt]{amsart}
\newtheorem{thm}{Theorem}[section]
\newtheorem{cor}[thm]{Corollary}
\newtheorem{lem}[thm]{Lemma}

\theoremstyle{definition}

\newcommand{\QED}{\rlap{$\sqcup$}$\sqcap$\smallskip}

\theoremstyle{remark}

\newcommand{\diam}{\operatorname{diam}}
\newcommand{\dist}{\operatorname{dist}}
\newcommand{\id}{\operatorname{id}}
\newcommand{\cl}{\operatorname{cl}}

\newcommand{\tl}{\tilde}
\newcommand{\eps}{\epsilon}
\newcommand{\di}{\partial}

\newcommand{\inter}{\operatorname{int}}
\newcommand{\orb}{\operatorname{orb}}

\newcommand{\ra}{\rightarrow}
\newcommand{\sm}{\setminus}
\newcommand{\Dil}{\operatorname{Dil}}
\newcommand{\HD}{\operatorname{HD}}

\newcommand{\BB}{{\cal B}}

\newcommand{\FF}{{\cal F}}
\newcommand{\HH}{{\cal H}}
\newcommand{\II}{{\cal I}}

\newcommand{\LL}{{\cal S\cal L}}

\newcommand{\QQ}{{\cal Q}}
\newcommand{\RR}{{\cal R}}

\newcommand{\YY}{{\cal Y}}
\newcommand{\SL}{{\cal S\cal L}}
\newcommand{\TT}{{\cal T}}
\newcommand{\Com}{{\cal C}\!{\it om}}
\newcommand{\Top}{{\cal T}\!\!{\it op}}
\newcommand{\QC}{{\cal Q}{\it C}}

\numberwithin{equation}{section}
\newsymbol\Subset 1362

\newcommand{\secref}[1]{\S\ref{#1}}
\newcommand{\lemref}[1]{Lemma~\ref{#1}}
\newcommand{\corref}[1]{Corollary~\ref{#1}}

\newcommand{\A}{{\Bbb A}}
\newcommand{\B}{{\Bbb B}}
\newcommand{\C}{{\Bbb C}}
\newcommand{\D}{{\Bbb D}}
\newcommand{\bH}{{\Bbb H}}
\newcommand{\Q}{{\Bbb Q}}
\newcommand{\R}{{\Bbb R}}
\newcommand{\J}{{\Bbb J}}

\newcommand{\T}{{\Bbb T}}
\newcommand{\U}{{\Bbb U}}

\newcommand{\comm}[1]{}

\def\IMSmarkvadjust{0 pt}
\def\IMSmarkhadjust{0 pt}
\def\IMSmarkhpadding{0 pt}
\def\SBIMSMark#1#2#3{
 \font\SBF=cmss10 at 10 true pt
 \font\SBI=cmssi10 at 10 true pt
 \setbox0=\hbox{\SBF \hbox to \IMSmarkhpadding{\relax}
                Stony Brook IMS Preprint \##1}
 \setbox2=\hbox to \wd0{\hfil \SBI #2}
 \setbox4=\hbox to \wd0{\hfil \SBI #3}
 \setbox6=\hbox to \wd0{\hss
             \vbox{\hsize=\wd0 \parskip=0pt \baselineskip=10 true pt
                   \copy0 \break%
                   \copy2 \break% 
                   \copy4 \break}}
 \dimen0=\ht6   \advance\dimen0 by \vsize \advance\dimen0 by 8 true pt
                \advance\dimen0 by -\pagetotal
	        \advance\dimen0 by \IMSmarkvadjust
 \dimen2=\hsize \advance\dimen2 by .25 true in
	        \advance\dimen2 by \IMSmarkhadjust

  \openin2=publishd.tex
  \ifeof2\setbox0=\hbox to 0pt{}
  \else 
     \setbox0=\hbox to 3.1 true in{
                \vbox to \ht6{\hsize=3 true in \parskip=0pt  \noindent  
                {\SBI Published in modified form as part of}\hfil\break
``Dynamics of quadratic polynomials, I-II'',
{\it  Acta Math.}~{\bf 178} (1997), 185-297.
 
                \vfill}}
  \fi
  \closein2
  \ht0=0pt \dp0=0pt
 \ht6=0pt \dp6=0pt
 \setbox8=\vbox to \dimen0{\vfill \hbox to \dimen2{\copy0 \hss \copy6}}
 \ht8=0pt \dp8=0pt \wd8=0pt
 \copy8
 \message{*** Stony Brook IMS Preprint #1, #2. #3 ***}
}

\begin{document}

\title[rigidity]{Dynamics of quadratic polynomials.\\
 II. Rigidity.}
\author {Mikhail Lyubich }
%\thanks{
%   This work was supported in part by Sloan Research Fellowship
% and NSF grants DMS-8920768 and DMS-9022140.}
%\date{October 22, 1995}

\maketitle

\thispagestyle{empty}
\SBIMSMark{1995/14}{October 1995}{}

\input{psfig}
\section{Introduction}\label{intro}

This is a continuation of the series of notes on the dynamics
of quadratic polynomials.  
The first part of this series \cite{L6} will be
 systematically used for the reference. In particular
we will assume that the reader is familiar with the background 
outlined in \S 2 of Part I: quadratic-like maps, straightening,
combinatorial classes, external rays, 
the Mandelbrot set $M$, secondary limbs, puzzle, etc.

Let $f$ be a quadratic-like map which does not have non-repelling
periodic points.
Let us say that $f$ satisfies the 
{\it secondary limbs condition} if there is a finite family of
truncated secondary limbs  $L_i$ of the Mandelbrot set such that
the hybrid classes of all renormalizations $R^m f$ belong to
$\cup L_i$. Let 
$\SL$ stand for the class of quadratic-like maps satisfying the
secondary limbs condition. 

Here are some examples of maps of class $\SL$:

\smallskip\noindent$\bullet$ Maps which are at most finitely 
renormalizable and don't have non-repelling periodic points
 (Yoccoz class); 

\smallskip\noindent$\bullet$ Infinitely renormalizable maps of  bounded
type (``bounded type" means that the relative periods of all renormalizations
are  bounded);

\smallskip\noindent$\bullet$ Real maps which don't have non-repelling periodic
points. 

\smallskip
Recall that a quadratic-like map $f$ has {\it a priori bounds} 
if there is an $\eps>0$ such that $\mod(R^m f)\geq \eps>0$
for all renormalizations.  

The goal of this paper is to prove the following result:
 
\proclaim Rigidity Theorem.  Any combinatorial class contains at
most one quadratic polynomial satisfying the secondary limbs 
condition with {\it a priori bounds}.

We believe that the second assumption actually follows from the first one:

\proclaim Conjecture. The secondary  limbs condition 
implies {\it a priori} bounds. 

This conjecture is supported by a few partial results (see below).
Note, however, that {\it a priori } bounds don't hold
for all quadratics: see examples of non-locally connected Julia sets
\cite{M1}.

Let $\QC(c)\subset\Top(c)\subset\Com(c)\subset \C$ 
stand  respectively 
for the quasi-conformal, combinatorial and topological classes of the
quadratic map $P_c$. A map $P_c$ is called combinatorially 
(respectively topologically or quasi-conformally) rigid if
$\Com(c)= \{c\}$ (respectively $\Top(c)=\{c\}$ or $\QC(c)=\{c\}$).
%Clearly, combinatorial rigidity implies topological rigidity,
%which in turn implies quasi-conformal rigidity.

\proclaim Corollary I.  Assume that all maps
of $\Com(c)$ (respectively $\Top(c)$) satisfy the secondary limbs
condition with
{\it a priori} bounds. Then  $P_c$ is combinatorially (respectively
topologically) rigid.

The corresponding quasi-conformal rigidity problem is settled
by  McMullen's Rigidity Theorem \cite{McM2}  which asserts that any
quadratic polynomial with {\it a priori} bounds is quasi-conformally
rigid.

The strongest, combinatorial,
 rigidity  of a map $P_c$ turns out to be equivalent
to the local connectivity of the Mandelbrot set $M$ at $c$ (see \cite{DH1,Sch}).
This property of $M$ was conjectured by Douady and Hubbard under the name 
``MLC". Prior to this work it was established in the  following cases:

\smallskip\noindent $\bullet$ Parabolic points (Douady and Hubbard \cite{DH1});

\smallskip\noindent $\bullet$  Boundaries of the hyperbolic components of $M$
   (Yoccoz, see Hubbard \cite{H});

\smallskip\noindent $\bullet$  At most finitely renormalizable maps
  (Yoccoz, see Hubbard \cite{H}, Kahn \cite{K}). \smallskip

The following  Corollary  adds a pool of infinitely renormalizable 
maps to this list.
In Part I of this paper {\it a priori bounds} have been proven 
for all maps of class $\SL$ with sufficiently big type 
(in the sense of Theorems IV and IV$'$ of Part I). Thus we have:

\proclaim Corollary II \cite{L6}.  A quadratic polynomial
$P_c\in \SL$ of a sufficiently big  type is rigid, so that
the Mandelbrot set is locally connected at $c$.

In particular, this gives first examples of infinitely
renormalizable parameter values $c\in M$ of {\it bounded type}
 where MLC holds (though one does not need the full capacity of Corollary II
to produce some examples of such kind). 

\smallskip {\it Remark.}  It is easy to construct some 
infinitely renormalizable parameter values of unbounded type
where MLC holds (oral communication by A. Douady). 
First find arbitrary small copies  $M_n$ of the
Mandelbrot set near $c=-2$.  Then for an appropriate subsequence
$n(k)$, the tuned  Mandelbrot copies
 $M_{n(1)} * M_{n(2)} *\dots * M_{n(l)})$ shrink to a single point.
\smallskip

One might wonder  of how big is the set of infinitely
renormalizable  parameter values
 satisfying the assumptions of Corollary II. We can show that  this set has
Lebesgue measure zero and 
Hausdorff dimension at least 1 (in preparation). Note
that 1=(1/2)2 where 
  $2=\HD(\di M)$ by Shishikura's Theorem \cite{Sh}.

Let us now dwell on the case of real parameter values. In this case,
"sufficiently big type" means sufficiently big {\it essential period}
(see \cite{LY} for the precise definition). For maps with
''small" essential period, the MLC problem is still open.
However, {\it a priori} bounds have been established for
all infinitely renormalizable real quadratics 
 (see \cite{S,MvS,L5,GS,LS,LY}).
Let us say that a parameter value  $c\in \R$ (or the corresponding
quadratic polynomial $P_c$) is {\it rigid on the  real line}
if $\Com(c)\cap\R=\{c\}$.
Thus we have:

\proclaim Corollary III. Any quadratic polynomial $P_c$ without
attracting cycles is rigid on the real line.

By the Milnor-Thurston kneading theory \cite{MT}, Corollary III implies:

\proclaim Corollary IV. Hyperbolic quadratics are dense on the real
line.

The last two Corollaries were first announced by Swiatek \cite{Sw}
who  approached them by methods of real dynamics. 
The methods of holomorphic dynamics presented in this paper
were developed in \cite{L5}.  

Another application of the above Rigidity Theorem is a construction
of the unstable manifolds for the renormalization operator at
infinitely renormalizable points of bounded type (in preparation). 

Let us now outline the structure of this paper. In \secref{J-l-c}
we show that the secondary limbs condition and {\it a priori}
bounds yields a definite space between the bouquets of little
Julia sets. This  provides us with special disjoint  neighborhoods
 of little Julia bouquets with bounded geometry (called ``standard").
  Also, together with the work of Hu \& Jiang \cite{HJ,J}
and McMullen \cite{McM3} this yields local connectivity of the
corresponding Julia set (Theorem I).
%This section also contains some preparation lemmas and constructions.

We start \secref{rigidity} with a  discussion of  reductions which
boil   the Rigidity Theorem down to the following problem:
Two topologically equivalent maps (satisfying the assumptions of
the theorem) are  Thurston equivalent. Then we set up an inductive
construction of  a sequence of approximations to the Thurston conjugacy. 
 In particular, we adjust   an approximate
conjugacy in such a way that it respects the standard neighborhoods
of  little Julia bouquets.   

The  last section, \secref{sec: principal nest},
 which  follows \S 4 of \cite{L5},
presents the proof of the Main Lemma. This  lemma gives a uniform
bound on the Techm\"uller distance between the
generalized renormalizations of two combinatorially equivalent
quadratic-like maps (the bound depends only on the selected secondary
limbs and {\it a priori} bounds). The main geometric ingredient
which makes this work is the linear growth of the principal moduli
proved in Part I of this paper.

In the Appendix we collect necessary background material in the theory of
quasi-conformal maps.

In conclusion let us make a couple of remarks on history and some 
related results and methods.
The origin of our approach to the rigidity problem can be tracked back to 
the proof of  Mostow Rigidity: from topological to quasi-conformal equivalence,
and then (by means of ergodic theory)  from quasi-conformal to conformal
equivalence. This set of ideas were brought to the iteration theory
by Sullivan and Thurston. 

The passage from quasi-conformal to conformal equivalence in our
setting  is settled by McMullen's Rigidity Theorem \cite{McM3}.
Our main task is to pass from topological to quasi-conformal equivalence.
A way to do this called ``pull-back argument" is to start with a 
quasi-conformal map respecting some dynamical data, and to pull it back 
so that it will respect more and more data on every step. In the
end it will become (with some luck) a quasi-conformal conjugacy.
This method was introduced by Thurston 
(see \cite{DH3} and also \cite{McM1})
for postcritically finite maps, and exploited by
Sullivan  \cite{S,MvS} for real infinitely renormalizable maps
of bounded type. These first applications dealt with maps with rather
simple combinatorics.

For  more complicated combinatorics, a certain real version of this method
based on the so called  ``inducing"  was suggested
by Jacobson \& Swiatek \cite{JS,Sw}. (Roughly speaking, ``inducing"
means building out of $f$ an expanding map with a definite range.)
% but it seems to be very hard to 
%carry it out without complex extension of the map and geometric
%preparation in the complex plane. 
On the other hand, 
by means of a purely
complex pull-back argument in the puzzle framework,
Jeremy Kahn \cite{K}  proved removability of non-renormalizable
Julia sets (which  yields the Yoccoz Rigidity Theorem) . 

Our way is  different from all the above, though it has some common features
with them.
We believe that holomorphic dynamics is the right framework for the
rigidity  problem, and our method is  purely complex.
Rather than building an induced expanding map,
 we pass consecutively from bigger to smaller scales by means of 
the generalized  renormalization \cite{L1}, and carry out the pull-back 
using growth of moduli and complex {\it a priori} bounds \cite{L5,L6}.  

Let us note that there is a different approach to  rigidity  problems,
by comparing the dynamical and parameter planes. This method was used
by Branner \& Hubbard \cite{BH} to prove rigidity of  cubic maps
with one escaping critical point and
``non-periodic tableaux" (which corresponds to non-renormalizable 
quadratics). It was also used by Yoccoz to prove  rigidity of 
at most  finitely renormalizable quadratics. 
In the forthcoming notes we will discuss this approach in our setting.

Let us also note that the MLC problem is closely related to 
the problem of landing of parameter rays at points $c\in \di M$.
MLC certainly yields landing of all rays, but, on the other hand,
landing of some special rays has  been a basis for progress in the 
MLC problem.
 The first results in
this direction (landing at parabolic and Misiurewicz points) were
obtained by Douady \& Hubbard (see \cite{DH1,M2,Sch}). 
Recently Anthony Manning \cite{Ma}  has estimated  
the Hausdorff dimension of  the set of rays landing
at infinitely renormalizable points.

\medskip {\bf Notations and terminology.} 
 Throughout the paper $f$ will stand for a 
quadratic-like map with critical point at 0.

Saying that a modulus of some annulus $A$ is {\it definite} means
that $\mod A\geq \eps>0$, where $\eps$ depends only on the selected 
truncated secondary limbs and
{\it a priori} bounds. Saying that some quantity is  {\it bounded} has 
an analogous meaning. 

Given a family of compact subsets $X_i\subset U$, we say that
 there is {\it a definite space} (at least $\mu>0$)
 in between them (in a domain $U$)
if for any $i$, there exists an annulus
 $A_i\subset U\sm \cup X_i$ with a definite modulus (at least $\mu$)
 which goes around
$X_i$ but does not go around other sets $X_j$, $j\not=i$. If $U$ is not
specified, then $U=\C$.

We will use the following notations:

\smallskip\noindent $\D_r=\{z: |z|<r\}$ is the standard disk of radius $r$, 
     $\D\equiv\D_1$ is the unit disk;

\noindent
$\T_r=\di \D_r$ is the standard  circle of radius $r$, $\T\equiv \T_1$ is the
    unit circle;

\noindent $\A(r,R)=\{z: r<|z|<R\}$ is a standard annulus;
  similar notation is used for a closed annulus
$\A[r,R]$ (or a  semi-closed one).

\noindent Let $P_c: z\mapsto z^2+c$.

\noindent As usual, $\omega(z)\equiv\omega(f,z)$ stands for the limit set
of the forward orbit $\{f^n z\}_{n=0}^\infty$. The set $\omega(0)$ is called
{\it postcritical}. 

\noindent  $R^m f$ is the $m$-fold renormalization of $f$.

\comm{
Let $J^m_i$ denote the little Julia sets of level $m$, that is,
 $J^m\equiv J^m_0=J(R^m f)$ 
and $J^m_i=f^i J^m$, $i=0,\dots , p_m-1$. 
They are organized in the pairwise disjoint bouquets
$B^m_j=B^m_j(f)$ of the Julia sets touching  at the same periodic  point.
Namely, if level $m-1$ is immediately renormalizable with period $l$ then
each $B^m_j$ consists of $l$ little Julia sets $J^m_i$ 
touching at their $\beta$-fixed points. 
Otherwise the bouquets $B^m_j$ just 
coincide with the  little Julia sets $J^m_j$.
  By $B^m\equiv B^m_0$ we will denote
the {\it critical} bouquet containing the critical point 0.
Let  $\J^m=\J^m(f)=\cup_i J^m_i=\cup_j B^m_j$.}

\medskip {\bf Acknowledgement.} I would like to thank Curt McMullen and
Yair Minsky for useful discussions.  I also thank MSRI for their hospitality:
Part of this work was done during the Complex Dynamics 
and Hyperbolic Geometry spring program 1995.
 This work has been partially supported
by NSF grants  DMS-8920768 and DMS-9022140, and the Sloan Research Fellowship. 

\bigskip
\goodbreak

\section{  Space between Julia bouquets.}\label{J-l-c}

\subsection{ Space and unbranching}\label{sec:space}
%Let us say that an annulus $A^m_j$ is associated with some $B^m_j$ if
%$A^m_j\subset \C\setminus \J^m$ and $A^m_j$
%goes around $B^m_j$ but does  not go around $B^m_l$ for $l\not= j$.  
%We say that {\it there is a definite space } in between the bouquets  if
%for any bouquet $B^m_j\subset \C\setminus  \J_m$ 
%there is an  annulus  $A^m_j$ with a definite modulus 
%which goes around this bouquet $B^m_j$ but does  not go around
%any other bouquet  $B^m_l$,  $l\not= j$. 
%%($\mod(A^m_j)\geq \mu>0$ is independent of $n,j$).
%
\noindent 
Let $J^m_i$ denote the little Julia sets of level $m$, that is,
 $J^m\equiv J^m_0=J(R^m f)$ 
and $J^m_i=f^i J^m$, $i=0,\dots , r_m-1$. 
They are organized in the pairwise disjoint bouquets
$B^m_j=B^m_j(f)$ of the Julia sets touching  at the same periodic  point.
Namely, if level $m-1$ is immediately renormalizable with period $l$ then
each $B^m_j$ consists of $l$ little Julia sets $J^m_i$ 
touching at their $\beta$-fixed points. 
Otherwise the bouquets $B^m_j$ just 
coincide with the  little Julia sets $J^m_j$.
  By $B^m\equiv B^m_0$ we will denote
the {\it critical} bouquet containing the critical point 0.
Let  $\J^m=\J^m(f)=\cup_i J^m_i=\cup_j B^m_j$.  Finally let $K^m_i$ 
be little filled Julia sets.

We will use the notation $F_m$ for the quadratic-like map $f^{p_m}$ near
any little Julia set $J^m_i$ (it should be clear from the context which
one is considered). In particular,
 $F_m=R^m f$ near the critical Julia set $J^m\ni 0$. 

Recall that $\QQ(\mu)$ stands for the space of quadratic-like maps
$f$  with  $\mod(f)\geq \mu>0$ supplied with the Caratheodory topology
(see \cite{McM2}  and  \S 5.6 of Part I). 
 Take a little copy  $ M'\subset M$  of the Mandelbrot set with root
at $b$.
Let   $\QQ(\mu , M')$  denote the subspace of $\QQ(\mu)$  consisting
of renormalizable quadratic-like maps $f$ 
such that the hybrid class of $Rf $ belongs to $M'\sm \{b\}$.
 
Let us have a family $\FF$ of  sets $X_a\subset  \C$ depending on some
parameter $a$ ranging over a topological space $\TT$. This dependence 
is said to be  (sequencially) 
{\it upper semi-continuous} if for any $a(i)\to a$,
the Hausdorff limit  of $X_{a(i)}$ is contained in $X_a$. For example it
is easy to see that the filled Julia set $K(f)$  of 
a quadratic-like map $f$ depends upper
semi-continuously on $f$. 
Let us say that a family $\FF$ of  sets $X_f\subset  \C$ is 
(upper) {\it semi-compact}
if any sequence $X_n$ of these sets contains a subsequence 
$X_{n(i)}$ converging in Hausdorff topology to a subset of some $X\in \FF$.

\begin{lem}\label{Caratheodory} 
The little filled Julia sets $K^1_i(f)$  form a semi-compact family
of sets as $f$ ranges over the space $\QQ( \mu, M')$.
\end{lem}

\begin{pf} By the Compactness Lemma (see \S 5.6 of Part I),
the space $\QQ(\mu, M')$ is compact.
Moreover the quadratic-like map $F_1$ depends continuously on 
$f\in \QQ(\mu, M')$ near any $K^1_i$. In turn, 
 the little  filled Julia sets $K^1_i$ 
   depend upper semi-continuously on $F_1$.
\end{pf}

\begin{lem}\label{space}
 Let $f$ be a quadratic-like map of class $\LL$ with complex {\it a priori}
  bounds. Then there is a definite space in between its bouquets $B^m_j$.
%depending only on the selected secondary limbs and {\it a priori} bounds.
\end{lem}

\noindent{\bf Proof.} %Let  $\mod(R^m f)\geq\eps>0,\; m=0,1\dots$.
Let us take a bouquet $B^m$. Let $\II^m$ stand for the set of indices
$j$ such that $B^{m+1}_j\subset B^m$. We will show first that 
there is a definite  annulus 
$$T^m\subset B^m\sm \bigcup_{j\in\II^m} B^{m+1}_j,$$
 which goes around $B^{m+1}$ but  does not go around  other
bouquets  $B^{m+1}_j$, $j\in \II^m$.

%Otherwise  it is nothing to prove as there is only one bouquet $B^{m+1}$ inside% $B^m$.
%Since $\mod (R^n f)\geq\alpha$ and the hybrid class of $R^n f$ is chosen from
%a truncated secondary limb, the initial Yoccoz puzzle pieces have a definite
%geometry. Moreover, if $R^n f$ is not immediately renormalizable then
%by \cite{L}, there is a definite space in between the little Julia sets
%$J^{n+1}_i$ of the next level contained in $J^n$. 
%Otherwise the next renormalization type is selected from only finitely many
%little Mandelbrot copies (corresponding to the choice
%of secondary limbs). Such family of quadratic-like maps with
%{\it a priori } bounds is compact. It follows that
%there is a definite space between the bouquets $B^{n+1}_j$ contained
%in $J^n$.

If $R^m f$ is not immediately renormalizable, then this follows from
 point (ii) of Theorem II  (Part I). 
So assume that $R^m f$ is immediately renormalizable. 

If $B^m=J^m$, then
it is nothing to prove as there is only one bouquet $B^{m+1}$ inside $B^m$.
Otherwise there are only finitely many renormalization types
producing the bouquet $B^m$ 
(which correspond to the little
Mandelbrot sets attached to the main cardioid and belonging to the selected
secondary limbs). 
%As $R^{m-1} f$ ranges over renormalizable
%quadratic-like maps with $\mod R^{m-1} f\geq \delta>0$,
% the bouquet $B^m$ ranges over a compact family of sets
%(in the Hausdorff metric).
By \lemref{Caratheodory},
 the bouquets $B^{m+1}_j$ contained in $B^m$ belong to a
compact family of sets. 
As they don't touch each other,
there is a definite space in between them. 

%Let us now show that there is a definite annulus  going around $B^{m+1}$
%which does not intersect and  does not go around any other 
% bouquets $B^{m+1}_j$, $j\not=0$.
Let $N(L, \eps)$ denote an $\eps \cdot\diam L$-neighborhood of a set $L$
(that is, the set of points on distance at most $\eps \diam L$ from $L$).
%Take a bouquet $B^{m+1}_j\subset B^{m}_s$. 
We have shown 
that there is an $\eps>0$ such that the neighborhood $N(B^{m+1},\eps)$
 does not intersect  other bouquets $B^{m+1}_j$ contained in the
same $B^{m}$. In particular, $N(B^1,\eps)$ does not intersect any
other $B^1_j$ (as all of them are contained in $B^0\equiv J(f)$). 
%This allows us to start the induction. 

Let us show by induction that 
\begin{equation}\label{empty intersection}
N(B^m,\eps)\cap B^m_k=\emptyset,\; k\not=0
\end{equation}
%Let us now consider a bouquet $B^{m+1}_k\not\subset B^{m}$.
Assuming this for  $m$, we should show that   
\begin{equation}\label{pass to m+1}
N(B^{m+1},\eps)\cap B^{m+1}_j=\emptyset,\; j\not=0.
\end{equation}
As we already know (\ref{pass to m+1}) for $j\in\II^m$, let $j\not\in\II^m$.
Then $B^{m+1}_j\subset B^m_k$ for some $k\not=0$,
 while $N(B^{m+1},\eps)\subset N(B^m,\eps)$,
and (\ref{pass to m+1}) follows from (\ref{empty intersection}).

What is left, is to show that there a definite space around 
any bouquet $B^{m+1}_j$
(not only around the critical one).
 But there is an iterate $f^l$ which univalently
maps $B^{m+1}_j$ onto $B^{m+1}$. 
 Pulling back the space around $B^{m+1}$ we obtain
the desired space about $B^{m+1}_j$.
\QED

An infinitely renormalizable map $f$ is said to
satisfy an {\it unbranched a priori bounds} condition (see \cite{McM3})
if for infinitely many levels $m$, there is a definite space in between
$J^m$ and the rest of the postcritical set, $\omega(0)\sm J^m$.

\begin{lem}\label{unbranching} A map $f\in\LL$ with {a priori} bounds
satisfies an unbranched a priori bounds condition.
\end{lem}

%\begin{lem}\label{space-1} For a map $f\in\LL$ with {a priori} bounds,
%there are infinitely many levels $m$ with a definite space in between
%$J^m$ and $\omega(0)\sm J^m$.
%\end{lem}

\noindent{\bf Proof.} We will show that the unbranched condition can fail
only if  the level $m$ is not immediately renormalizable, while
  $m-1$ is immediately renormalizable.
 As the complimentary sequence of levels
is infinite, the lemma will follow.

 If $R^{m-1}f$ is not immediately renormalizable
then the bouquet $B^m$ coincides with the little Julia set $J^m$.
By \lemref{space}, there is a definite space in between $J^m$ and 
$\J^m\sm J^m$. As $\omega(0)\sm J^m\subset\J^m\sm J^m$, the unbranched
condition holds on level $m$.

Assume now that both levels $m-1$ and $m$ are immediately
renormalizable. Then we will show that there is a definite space
in between $J^m$ and $\BB^{m+1}\equiv \cup_{j\not=0} B^{m+1}_j$. 

By \lemref{space}, there is a definite space in between $B^m\supset J^m$ and
$\BB^{m+1}\sm B^m$. So we should  check that there is
 a definite space in between $J^m$ and  $\BB^{m+1}\cap B^m$
(that is, the union of
non-critical bouquets $B^{m+1}_j$ contained in $B^m$). 
But $J^m$ does not touch any such $B^{m+1}_j$. Indeed, the only point
where they can touch could be the $\beta$-fixed point $\beta_m$ of $J^m$. 
But one can easily see that the little Julia sets of level
 $m+1$ never contain $\beta_m$.
By \lemref{Caratheodory}
there is a desired  space. 

Finally, as $\omega(0)\sm J^m\subset
\BB^{m+1}$, the statement follows.  \QED

{\bf Remark.} If  $R^m f$ is
not immediately renormalizable, while
$R^{m-1}f$ is immediately renormalizable, then  
the unbranched condition can fail. Indeed in this case there are 
several Julia sets $J^m_i$ which touch at the common fixed point
$\beta_m\in J^m$.  But the postcritical set $\omega(0)\cap J^m_i$
  can come arbitrarily close to $\beta_m$ 
(when $R^{m} f$ is a small perturbation of
a map  whose critical orbit eventually lands at $\beta_m$).

\subsection{ Local connectivity of Julia sets.}\label{JLS-sec}
Hu and Jiang  [HJ] proved that the Feigenbaum quadratic polynomial has
locally connected Julia set. The proof makes use of Sullivan's {\it a priori}
bounds (see \cite{MvS,S}). 
Then a more general result of this kind was worked out:
Any infinitely renormalizable quadratic map with unbranched 
{\it a priori} bounds
 has locally connected Julia set (see \cite{J,McM3}).
 Together with \lemref{unbranching} this yields (compare Theorem V of Part I):

\proclaim Theorem I. Let $f\in \LL$ be an infinitely renormalizable quadratic
polynomial with {\it a priori } bounds. Then the Julia set $J(f)$ is locally connected.
In particular, all maps from Theorems IV and IV$'$ of Part I have locally connected
Julia sets. 

For the sake of completeness, we will give a proof of this result.

\noindent{\bf Proof.}
{\it A priori } bounds imply that the ``little" Julia sets $J^m$ 
shrink down to the critical point. Indeed let 
$f_m\equiv R^m f\equiv f^{p_m}: U_m'\rightarrow U_m$ where
 $\mod(U_m\backslash U_m')\geq\eps>0$, with an $\eps$ independent of $m$.
Clearly $U_m$ does not cover the whole Julia set.

Let $\Gamma_m\subset U_m\backslash U_m'$ be a horizontal curve in the
annulus $U_m\backslash U_m'$ which divides it into two sub-annuli of
modulus at least $\epsilon/2$, and $\Gamma_m'\subset U_m'$ be its pull-back
by $f_m$.  By the Koebe Theorem, these curves have a 
bounded eccentricity about 0 (with a bound depending on $\epsilon$).
Since the inner
radius of  curve $\Gamma_m'$ about $0$ tends to 0 as $m\to\infty$
(it follows from the fact that the sufficiently high iterates
 of any disk intersecting $J(f)$ cover the whole $J(f)$), 
the $\diam\Gamma'_m\to 0$
as well. All the more, $\diam(J_m)\to 0$ as $m\to\infty$.
 
Let us take a $\delta>0$, and find an $m$ such that $J_m$ is
contained in the $\D_\delta$.

Let us now inscribe into $\D_\delta$ a domain bounded by equipotentials and
external rays of the original map  $f$.
Let $\alpha_m$ denote
the dividing fixed point of the Julia set $J^m$, and $\alpha_m'=-\alpha_m$
be the symmetric point. Let us consider a puzzle piece $P^{m,0}\ni 0$
bounded by any equipotential and four
external rays of the {\it original map $f$}
 landing at $\alpha_m$ and $\alpha_{m}'$. This is 
a ``degenerate" domain of the renormalized map $F_m$ (see \S 2.5 of Part I).
 By definition of the renormalized Julia set, the preimages 
$P^{m,k}\equiv F_m^{-k} P^{m,0}$  shrink down to  $J^m$.
%Let us then pull these puzzle-pieces back in the usual way,
%and use the notation $P^{m,l}(a)$ for the puzzle-piece of level $l$
%containing a point $a$.
Hence there is a puzzle piece $P^{m,l}$
 contained in the $\D_\delta$. As $J(f)\cap P^{m,l}$ is clearly connected,
the Julia set $J(f)$ is locally connected at the critical point.

Let us  now prove local connectivity at any other point $z\in J(f)$.
This is done by a standard spreading  of the local information
near the critical point  around the whole dynamical plane. 
%As we have seen, for any unbranched level $m$, 
%there is a choice of the renormalized map 
%$f_m: U_m'\rightarrow U_m$, 
%such that $(U_m\backslash J_m)\cap\omega(0)=\emptyset$,
%$\mod(U_m\backslash U_m')\geq\eps>0$, with an $\eps$ independent of $m$,
%and the $U_m, U_m'$ have bounded geometry.
Let us consider two cases.\smallskip

{\sl Case (i)}. Let the orbit of $z$ accumulates on all Julia sets
  $J^m$. Let $m$ be an unbranched level.
Then there is  an $l=l(m)$ such that  the puzzle piece
 $P^{m,l}$ is well inside $\C\sm(\omega(0)\sm J^m)$. 

 Take now
the first moment $k=k(m)\geq 0$ such that $f^kz\in P^{m,l}$.
Let us consider the pull-backs $ Q^{m,l}\ni z$ 
of $ P^{m,l}$
 along the orbit $\orb_k(z)=\{z,...,f^k z\}$. 
By Lemma 3.3 of Part I, this pull-back is univalent.
Moreover, it allows a univalent extension to a definitely bigger domain.
%Thus $U_m\sm V_m'$ provides us with a definite Koebe space for
%the corresponding branch of $f^{-k}$.

 By the Koebe Theorem,  $Q^{m,l}$ has a bounded eccentricity about
$z$. Since the inner radius of this domain about $z$ tends to 0 as $m\to\infty$,
 the $\diam Q^{m,l}\to 0$ as well. 
As $Q^{m,l}\cap J(f)$ are connected, the Julia set is locally connected
at $z$.

\smallskip{\sl Case (ii)}. Assume now that the orbit of $z$
 does not accumulate on  some $J^m$. 
Hence it accumulates on some point $a\not\in \omega(0)$.
  Let us consider the puzzle associated with the periodic point $\alpha_m$
(so that the initial
 configuration consists of a certain  equipotential and the
 external rays landing at $\alpha_m$). 
 Since the critical puzzle pieces shrink to 
$J^m$,  the puzzle pieces  $Y^{l}_i$ of sufficiently big depth $l$
containing $a$ are
 disjoint from $\omega(0)$ (there are several such pieces if  $a$
is a preimage of $\alpha_m$). Take such an $l$, and 
let $X$ be the union of these
puzzle pieces. It is a closed topological disk disjoint from $\omega(0)$ whose
interior contains  $a$.  
%Hence there is a simply connected
%neighborhood $V\supset X$ still disjoint from
%$\omega(0)$.

 Consider now the moments $k_i\to\infty$ when the orbit of
$z$ lands at $\inter X$, and pull $ X$ back to $z$. 
By the same Koebe argument as in case (i) we conclude that these
 pull-backs shrink to $z$. 
It follows that $J(f)$ is locally
connected at $z$.           \QED

\subsection{ Standard neighborhoods.}\label{sec:standard neighborhoods}
In this section  we will construct some special fundamental domains
near little Julia bouquets. 
Let us consider first 
the non-immediately renormalizable case when the construction
can be done in a particularly nice geometric way. 

\begin{lem}\label{geometric annuli}
 Let $f$ be $m$ times renormalizable quadratic map. Assume that
 the space in between the little Julia sets $J^m_i$ is at least $\mu>0$.  Then there are
  disjoint fundamental annuli $A^m_i$ around little Julia sets $J^m_i$, with 
$\mod A^m_i\geq\nu(\mu)>0$.
\end{lem}

\noindent{\bf Proof.} Let us consider the Riemann surfaces 
$S=\C\sm\J^n$ and
  $S'=\C\sm f^{-1}\J^n\subset S$.
 Then $f: S'\rightarrow S$ is a double branched
covering. Let us uniformize $S$, 
that is represent it as the quotient $\HH^2/\Gamma$
of the hyperbolic plane modulo the action of a Fuchsian group. In this conformal
representation $S$ admits a compactification $  S\cup \partial S$ to a bordered
Riemann surface, with the components $\partial S^m_i$ 
of the ''ideal boundary" $\partial S$
corresponding to the little Julia sets $J^m_i$.

Let $\hat S=S\cup\partial S\cup \bar S$ be the double of $S$, that is
 $(\C\backslash \Lambda(\Gamma))/\Gamma$, where $\Lambda(\Gamma)\subset S^1$ is the limit set
of $\Gamma$. The boundary components $\partial S^m_i$ are geodesics in $\hat S$. Moreover,
these geodesics have  hyperbolic length bounded by a constant $L=L(\mu)$ independent of $m$.
% Indeed by Lemma 6.1,
% there is a  definite space in between the little Julia sets $J^m_i$, and hence there is
%a definite collar around $\partial S^m_i$ in $\hat S$.

Let $\sigma: S\rightarrow S$ be the natural anti-holomorphic involution of $S$.
Let  $\bar S'=\sigma S'$ and $\hat S'=S\cup\partial S\cup \bar S'\subset \hat S$ 
be the double of  $S'$ inside $S$. Then $f$ allows an extension to a holomorphic
double branched covering $\hat f: \hat S'\rightarrow \hat S$ 
commuting with the involution
$\sigma$. Its restrictions $\hat f|\partial S^m_i$ 
are the double branched coverings
of the topological circles $\partial S^m_i$.

Let $C^m_i(r)\supset \partial S^m_i $ 
stand for the hyperbolic $r$-neighborhood of the
geodesic $\partial S^m_i$.
By the Collar Lemma (see [Ab]), there is an $r=r(L)$ 
(independent of the particular
Riemann surface and geodesics) such that the collars
 $C^m_i\equiv C^m_i(r)$ are pairwise disjoint. 
Moreover,  $\mod(C^m_i)\geq\mu(L)>0$.
% and there is a definite space $\mu(L)$ in between them. 

Let us now take such a collar $C=C^m_i$,  and let $\Gamma=\partial S^m_i$.
Let $C'\subset S'\cap C$ be the component of
$\hat f^{-1} S$ containing $\Gamma$. 
Then $\hat f: C'\rightarrow C$ is a double covering
preserving $\Gamma$. As we have in the hyperbolic metric of $S$:
$$\int_{\Gamma}\|D\hat f\|=2l(\Gamma),$$
there is a point $z\in\Gamma$ such that $\|Df(z)\|\geq 2$. 
This easily implies that
$\|D\hat f^{-1}(\zeta)\|\leq q(a)<1$ 
if the hyperbolic distance between $z$ and $\zeta$ does not
exceed $a$. In particular, $\|D\hat f^{-1}\|(\zeta)\leq q=q(L,r)<1$
 for all $\zeta\in C$. 

It follows that $C'$ is contained in the hyperbolic 
$r/q$-neighborhood of $\Gamma$,
and hence the annulus $\mod(C\backslash C')\geq \rho(r,q)=\rho(\mu)$. Let now 
$A^m_i=(C\backslash C')\cap S.$ \QED

Note that in the above lemma we don't assume 
{\it a priori} bounds but just a definite
space between the Julia sets (which thus implies {\it a priori } bounds). 
%The price we paid is the assumption that $R^{m-1}$ is not
%immediately renormalizable. 
Assuming {\it a priori} bounds,
let us now give a different construction which works in
the immediately renormalizable case as well.

\begin{figure}[htp]
\centerline{\psfig{figure=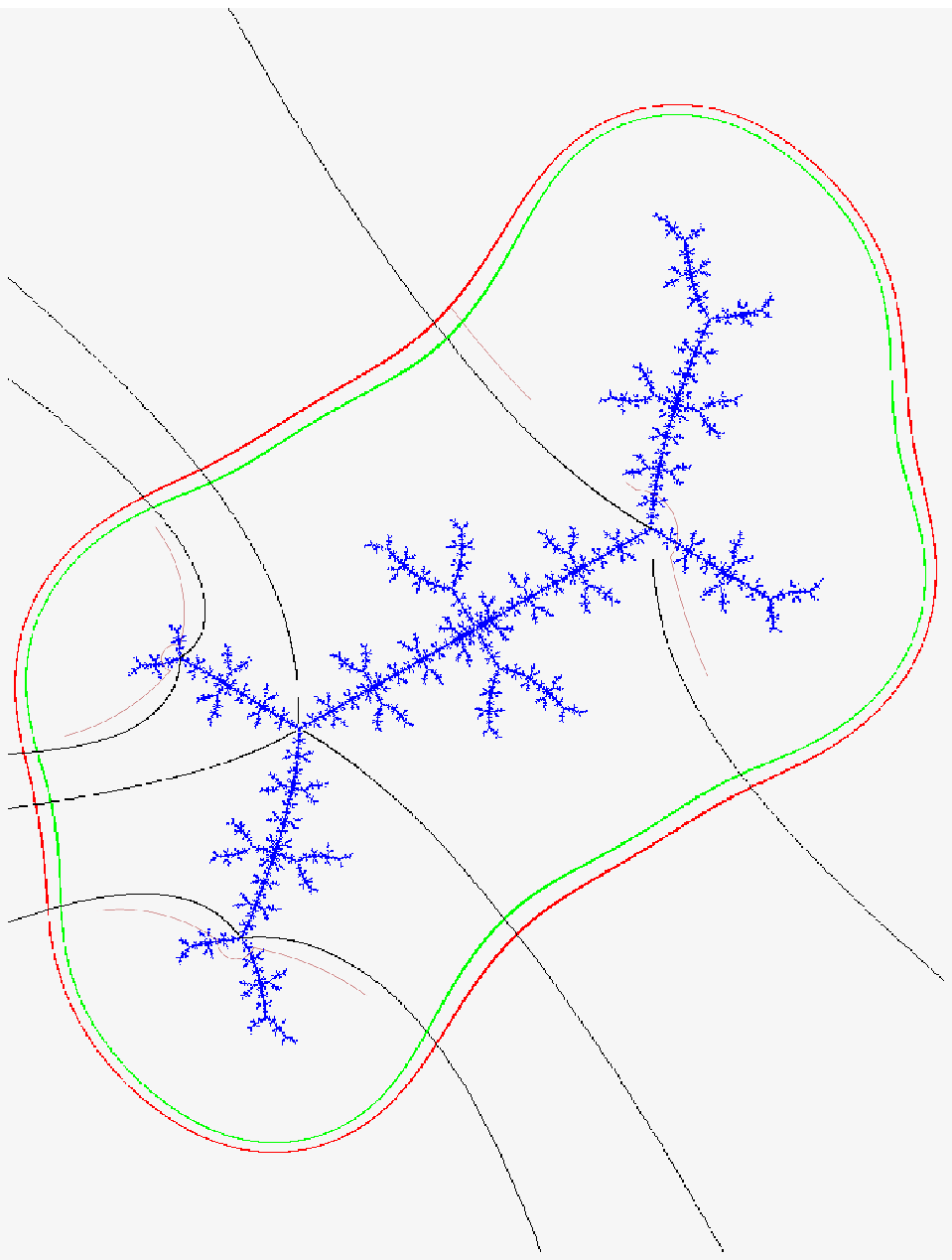,width=.8\hsize}}
\bigskip\centerline{Figure 1. Standard neighborhood of a Julia bouquet
 (made by B. Yarrington).}\bigskip
\end{figure}

Let us consider a bouquet $B^m_j=\cup_i J^m_i $ of level $m$, where $J^m_i$
 touch at point $\alpha_{m-1}$.
Let $b_{m,i}\in J^m_i$ be the points $F_m$-symmetric to  $\alpha_{m-1}$, 
that is,  $F_m b^m_i=\alpha_{m-1}$ (''co-fixed points"). 
 Let us consider the domain $\Upsilon^m_j$ bounded by the
pairs of rays  landing at these points
 (defined via a straightening of $F_{m-1}$), and $p_m$
arcs of equipotentials.
 Let us then thicken this domain near the points $b_i^m$ as
described in \S 2.5 of Part I 
(that is, replace the rays landing at $b^m_i$ by nearby rays and little
circle arcs around $b^m_i$). 
Denote the thickened domains by $U^m_j$ (see Figure 1).
We also require that
these domains are naturally  related by dynamics so that
$f \Upsilon^m_j=\Upsilon^m_k$ and
$f U^m_j=U^m_k$ whenever $f B^m_j=B^m_k$ and $B^m_j$ is non-critical.
 Let us call  $U^m_j$ a {\it standard neighborhood} of the bouquet $B^m_j$. 
Let $\U^m=\cup U^m_j$.

\begin{lem}\label{standard neighborhoods} 
Let $f$ be an $m$ times renormalizable quadratic map of class
$\SL$ with {\it a priori} bounds. Then 
 there exist disjoint standard neighborhoods
$U^m_j$ of $B^m_j$ with bounded geometry, and such that the annulus
$\mod(U^m_i\backslash B^m_i)$ have a definite modulus. 
\end{lem}

\noindent{\bf Proof.} By the Straightening Theorem,
the renormalization $R^{m-1}f: J^{m-1}_k\rightarrow J^{m-1}_k$ 
is $K$-qc conjugate
to a quadratic polynomial $P_c: z\mapsto z^2+c$, with $K$ dependent only on
{\it a priori bounds}. 
%Remember  that $c$ belongs to a finite family of truncated
%secondary limbs.

 Let $B\subset J(P_c)$ be the critical bouquet of little Julia sets of $RP_c$. 
Let $\Omega(\eps)$ be its  neighborhood 
 bounded by arcs of equipotentials of level
$1-\eps$, circle arcs  of radius $\eps$, 
and rays with arguments $\theta+t(\eps)$ (see \S 2.4 of Part I).
Here
$\theta$ are the arguments of the rays landing at the co-fixed points, and 
$t(\eps)\in (-\eps,\eps)$ is selected in such a way that $\Omega(\eps)$ 
is a renormalization
domain for any $P_c$ from selected truncated secondary limbs.

The geometry of these domains depends only on the selected limbs and $\eps$.
Also,  the Hausdorff distance $d_c(\eps)$ of $\di\Omega(\eps)$ to
$B$ tends to 0 as $\eps\to 0$ uniformly over 
$c$ belonging to the selected truncated
limbs. Indeed, this is clearly true for a given parameter value $c$. 
Take a little
$\delta>0$, and  find an $\eps=\eps_c$ such that $d_c(\eps_c)<\delta$. 
Then for all
$b$ sufficiently close to $c$, $d_b(\eps_c)<2\delta$. Compactness of
the truncated limbs completes the argument.    
 
It follows that for all  sufficiently small $\eps$
(depending only on the selected limbs and {\it a priori} bounds),
$\Omega(\eps)$ belongs to the range of the straightening map. Hence these neighborhoods
can be transferred to the dynamical
$f$-plane. We obtain  neighborhoods $U(\eps)$ of the corresponding bouquet $B$
with bounded geometry (depending on parameter $\eps$).

Moreover, as quasi-conformal maps are quasi-symmetric
(see Appendix),
 the Hausdorff distance from $\di U(\eps)$ to the bouquet $B$ is at most
$\rho(\eps)\cdot \diam B,$ where $\rho(\eps)\to 0$ as $\eps\to 0$.
 Hence for all sufficiently
small $\eps$, the neighborhood $U(\eps)$ is well inside the domain 
$\C\setminus \cup_{j\not=0} B^m_j$.

Let us now pull this neighborhood back by dynamics 
to obtain standard neighborhoods
$U^m_j(\eps)$ of other bouquets $B^m_j$. Since $U(\eps)$ is well inside  
$\C\setminus \cup_{j\not=0} B^m_j$, these pull-backs have a bounded distortion.
Hence the Hausdorff distance from $\di U^m_j(\eps)$ 
to the bouquet $B^m_j$ is at most
$\rho(\eps)\cdot\diam B,$ where $\rho(\eps)\to 0$ as $\eps\to 0$. 

Since by \lemref{space} there is
a definite space between the bouquets, there is also a definite space
between the  neighborhoods
$U^m_j(\eps)$,  for  for all $\eps\in (0,\eps_*]$ (with $\eps_*$
depending only on the selected limbs and {\it a priori} bounds). 
Also,  the moduli of $U^m_j(\eps)\setminus B^m_j$ depend only on the limbs,
{\it a priori} bounds and $\eps$. 
So they are definite, for instance in the range  
$\eps\in(0.01\eps_*, \eps_*]$. \QED

%\smallskip{\it Remark.} Let us take  neighborhoods
%$L^m_i\equiv U^m_i(\eps_*)\supset U^m_i(0.5\eps_*)\equiv U^m_i$. 
%Then clearly the
% annuli $L^m_i\setminus N^m_i$ and $U^m_i\setminus B^m_i$ 
%have a definite modulus.
%This remark will be useful.

We keep using the notations  $B^m_j$, $\Upsilon^m_j$ etc.
introduced before \lemref{standard neighborhoods},
and we also assume that the standard neighborhoods $U^m_j$ satisfy the
conclusions of \lemref{standard neighborhoods}.
We will define a special qc map
\begin{equation}\label{standard straightening}
S_{m}: (U^m_j\sm B^m_j)\ra \A(1,4). 
\end{equation}
with bounded dilatation. This map will be called a 
{\it standard straightening}, or a {\it standard local chart} near 
the bouquet $B^m_j$.

It follows from  {\it a priori bounds} assumption that
for any Julia set $J^l_i$ there exist
  Jordan disks $\Omega^l_i\supset\Pi^l_i\supset J^l_i$  
such that $F_l: \Pi^l_{i}\ra \Omega^l_i$ is a quadratic-like map, and there
exists a qc map
\begin{equation}\label{annuli straightening}
\Psi_{l,i}: (\Omega^l_i\sm J^l_i,\; \Pi^l_{i}\sm J^l_i)\ra (\A[1,4],\; \A[1,2])
\end{equation} 
 with bounded dilatation conjugating $F_l: \Pi^l_{i}\ra \Omega^l_i$ and
$P_0: \A[1,2]\ra \A[1,4]$, $P_0: z\mapsto z^2$.

 Moreover, if $J^m_i$ does not touch other 
Julia sets of level $m$  (that is, $F_{m-1}$ is not immediately 
renormalizable) then one can select the standard neighborhood
$U^m_i$ as $\Omega^m_i$. In this case let  us define the
standard straightening (\ref{annuli straightening}) as $\Psi_{m,i}$.

%Let us first assume that $F_{m-1}$ is not immediately renormalizable,
%so that $B^m_j$ consists of a single Julia set $J^m_j$.
%Then $U^m_j$ is a domain of the quadratic-like map
%$F_m$ such that the annulus  $A=F_m^{-1}U^m_j\sm U^m_j$ has a bounded
%geometry. Hence there is a qc map of $A$ onto  $\A[2,4]$ respecting 
%dynamics on the boundaries. Pulling this map back to the
%preimages of $A$, we obtain a desired map (\ref{standard straightening}).

If $F_{m-1}$ is immediately renormalizable, then let us consider a family 
of little Julia sets and bouquets:  
\begin{equation}\label{inclusions}
 \bigcup_i J^m_i= B^m_j\subset J^{m-1}_k.
\end{equation}
Let us cut   $\Upsilon^m_j$ by the rays landing at the fixed point 
$\alpha_{m-1}$ into components
 $\Xi_i^m\supset J^m_i.$ 
Since the hybrid class of $F_{m-1}$ may belong to a bounded number of 
little Mandelbrot sets (attached to the main cardioid and intersecting
the selected secondary limbs), the domains $\Xi^m_i$ have a bounded 
geometry. 
 Hence the maps $\Psi_{m,i}$ 
can be selected in such  a way that they have  bounded dilatation and
$$\Psi_{m,i}|\bigcup_i \Xi^{m}_i=\Psi_{m-1,k}.$$
Thus they glue together into  a single qc map  (\ref{standard straightening}).
%\begin{equation}\label{standard straightening}
%S_{m}: (U^m_j\sm B^m_j)\ra \A(1,4). 
%\end{equation}
%with bounded dilatation. 
%The map (\ref{standard straightening}) will be called a 
%{\it standard straightening}
%near the bouquet $B^m_j$, or a {\it standard local chart} near $B^m_j$.

By the rays and equipotentials near the bouquet we will mean the
 $S_m$-preimages
 of the vertical intervals and horizontal circles in the cylinder $\A(1,4)$.
This will be also referred to as the {\it standard coordinate system}
 near $B^m_j$.

Let us show in conclusion that the  little Julia bouquets and
corresponding standard neighborhoods exponentially decay.
Let $\diam (X)$ stand for the Euclidean diameter
of a set $X$.

\begin{lem}\label{exp decay} Let $f\in \SL$ be a quadratic-like map
with {\it a priori} bounds.
Then there exist constants   $\lambda<1$  and $l_0>0$ depending on the
choice of limbs and {\it a priori} bounds such that 
for any two Julia bouquets $B^{m+l}_j\subset B^m_i$,
$$\diam B^{m+l}_j\leq \lambda^l \diam J^m_i, \quad l\geq l_0$$
\end{lem}

\begin{pf} 
  Let us straighten the renormalization $R^m f$ near $J^m_i$
to a quadratic polynomial $P_c$. 
The dilatation  $K$ of the straightening depends only on  the
{\it a priori} bounds, and $K$-qc  maps are H\"older
continuous with exponent  $1/K$ (see \cite{A}).
 Hence it is enough to show that
for the quadratic map $P_c$, there exist
constants   $\lambda<1$  and $l_0>0$ depending on the
choice of limbs and {\it a priori} bounds such that

\begin{equation}\label{polynomial case}
\diam B^l_j\leq \lambda^l.
\end{equation}
(Now $B^m_j$, $\J^m$ etc.  stand for the objects associated to $P_c$).

Note that $J(P_c)\subset \D_2$. Let  $\rho_l$ be the hyperbolic metric
on $\D_3\sm \J^l$.
Let  $\gamma^m_i$ be the hyperbolic 
geodesic in $\D_3\sm \J^l$ homotopic to 
a curve $\Gamma^m_i\subset \C\sm \J^m$ going once around $B^m_i$ but not
going around other Julia bouquets $B^m_k$, $k\not=i$. 

By \lemref{space}, there are annuli $A^n_i\subset \D_3\sm \J^m$ in the homotopy
class of $\Gamma^m_i$ with
a definite modulus, $\mod(A^m_i)\geq\nu$. 
Let us pick $\Gamma^m_i$ as the hyperbolic geodesic in $A^m_i$.
Then the hyperbolic length of this geodesic in $A^m_i$ is at most
$\pi/\nu$. All the more, the hyperbolic length of $\gamma^m_i$
in $\D_3\sm \J^l$ is bounded by the same constant.

By the Collar Lemma (see \cite{Ab}), there are exist disjoint annuli
$\Lambda^m_i\subset \D_3\sm \J^l$ in the homotopy class of 
$\gamma^m_i$ with $\mod(\Lambda^m_i)\geq \eta=\eta(\nu)>0$. 
By the Gr\"otzsch inequality, $\mod (\D_3\sm B^l_j)\geq l\eta.$
Hence there is an absolute constant $C$ such that
$\diam B^l_j\leq Ce^{-l\eta}$  (see Appendix A1 in Part I), 
and  (\ref{polynomial case}) follows. 
\end{pf} 

\begin{cor}\label{exp decay of U} Under the assumptions of \lemref{exp decay},
there exist constants   $\lambda<1$  and $l_0>0$  such that 
for the standard neighborhoods $U^m_j\subset U^m_i$ the
following estimates holds:
$$\diam U^{m+l}_j\leq \lambda^l \diam U^m_i,\quad l\geq l_0.$$
\end{cor}

\begin{pf} Indeed, the standard neighborhoods $U^m_i$ are commensurable
with the corresponding Julia sets $J^m_i$.
\end{pf}    

\subsection{Removable sets} The reader is referred to the Appendix for the
definition and a discussion of removability.

\begin{lem}[McM2]\label{omega removable}
Under the assumptions of \lemref{exp decay}, the post-critical set 
$\omega(0)$ is a removable Cantor set coinciding with $\cap\J^m$.
\end{lem}

% Let $I^m_i=\omega(0)\cap J^m_i$. 
\begin{pf} 
It was shown in the proof of 
\lemref{exp decay} that for any $z\in  \omega(0)\subset \cap\J^m$,
there is a nest of disjoint annuli around $z$ with a definite modulus.
Thus the first statement follows from the Removability Condition (see Appendix). 

Clearly, $\omega(0)\subset \cap \J^m\subset \cap\U^m$.
Vice versa, by \lemref{exp decay},  $\cap \J^m$ is covered by the 
uniformly shrinking  bouquets $B^m_i$. 
As every $B^m_i$ contains a postcritical
point, $\omega(0)$ is dense in $\cap\J^m$.
 \end{pf}

Let us finish this section with stating a
   standard fact on removability of expanding Cantor sets.
Let $\{U_i\}$ be a finite family of closed topological disks  with disjoint
closures.
Let us consider a  Markov  map $g: \cup U_i\ra \C$ satisfying the following
property: If $\inter (g U_i\cap U_j)\not=\emptyset$ then $g U_i\supset U_j$.
 As usual, let 
$$K(g)=\{z:\; g^nz\in \cup U_i,\; n=0,1,\dots\}$$
stand for the filled Julia set of $g$.

\begin{lem}\label{Markov maps} For a Markov map as above,
the filled Julia set $K(g)$ is removable.
\end{lem}

\begin{pf} Let us select a family of annuli $A_j\subset g U_j\sm \cup U_i$
homotopic to $\di (g U_i)$ in  $g U_j\sm\cup U_i$.
Let consider cylinder sets
$U^m_{i(0), i(1),\dots, i(m-1)} $ defined by the following property:
$$g^k  U^m_{i(0), i(1),\dots, i(m-1)}\subset U_{i(k)},\; k=0,1,\dots, m-2;
 \quad g^{m-1} U^m_{i(0), i(1),\dots, i(m-1)}=U_{i(m-1)}.$$
The pull-back of the annulus $A_j$  to 
$U^m_{i(0), i(1),\dots, i(m-1)} \sm 
\bigcup_i U^{m+1}_{i(0), i(1),\dots, i(m-1), i}$
by the univalent map $g^m: U^m_{i(0), i(1),\dots, i(m-1)} \ra g U_{i(m-1)}$ 
has the same modulus as $A_{i(m-1)}$. 
This provides us with a nest of disjoint annuli with definite
modulus about  any $z\in K(g)$. The Removability Condition concludes the proof.
\end{pf}

\comm{
\begin{lem}\label{residual Cantor sets} For a quadratic polynomial
$f=P_c$ (perhaps, infinitely renormalizable), let us consider a 
critical puzzle piece $V\ni 0$ for some renormalization $R^m f$.
Then the set $$\Lambda_V=\{z:\; f^n z\not\in V,\; n=0,1,\dots\}$$
is removable.
\end{lem}

\begin{pf} The Cantor set $\Lambda\equiv\Lambda_V$ is expanding (see \cite{L..}
or \cite{Pr}). By the general theory of expanding maps (see \cite{Kz}),
$f$  is Markov in a neighborhood of $\Lambda$.
 By \corref{Markov maps},  $\Lambda$ is removable.
\end{pf}   }

\section{ Rigidity: beginning of the proof} 
\label{rigidity}

\subsection{Reductions} \label{statement}
In this section we begin  to prove the
 Rigidity Theorem stated in the Introduction.
  Since quadratic polynomials label hybrid classes of quadratic-like maps, 
this theorem can be stated in the following way:

\proclaim Rigidity Theorem (equivalent statement). 
Let $f, \tl f\in \LL$ be two quadratic-like maps with 
{\it a priori bounds}. If $f$ and $\tl f$ are combinatorially equivalent then 
they are hybrid equivalent.

The proof is split into three steps:

\smallskip\noindent Step 1. $f$ and $\tl f$ are topologically equivalent;

\smallskip\noindent Step 2. $f$ and $\tl f$ are qc equivalent;

\smallskip\noindent Step 3. $f$ and $\tl f$ are hybrid equivalent.\smallskip

The first step (passage from combinatorial to topological equivalence)
follows from the local
connectivity of the Julia sets (Theorem I).  Indeed, a locally connected Julia set is homeomorphic
to its combinatorial model (see \cite{D}). Since the combinatorial model is
the same over the combinatorial class, the conclusion follows.

The last step (passage from qc to hybrid equivalence)
is taken care of  McMullen's Rigidity Theorem \cite{McM2}.
Indeed, it asserts that an infinitely renormalizable quadratic-like 
map with {\it a priori} bounds 
does not have invariant line fields on the Julia set. 
It follows that if $h$ is a qc
conjugacy between $f$ and $\tl f$ then $\di\bar h=0$ 
almost everywhere on the Julia set.
Thus $h$  is a hybrid conjugacy between $f$ and $\tl f$.

So, our task is to take care of Step 2: 

\proclaim Theorem II.  Let $f, \tl f\in \LL$ be two quadratic-like maps
 with  {\it a priori
bounds}. If $f$ and $\tl f$ are topologically  equivalent then  
they are qc equivalent.

In what follows we will mark with tilde the objects for 
$\tl f$ corresponding to those
for $f$ (unless another meaning is explicitly assumed).
 When we introduce some objects for $f$, we  assume that the
corresponding tilde-objects are automatically introduced as well. 

\subsection{Thurston's equivalence}\label{Thurston}

Let $f: U'\ra U$ and $\tl f: \tl U'\ra \tl U$ 
 be two topologically equivalent quadratic-like maps.
% and $\psi$ be a conjugacy between them.
 Let us say that $f$ and $\tl f$ are {\it Thurston equivalent}
 if for  appropriate choice of
domains $U,U',\tl U, \tl U'$, there is a qc map
 $h:(U,U',\omega(0))\rightarrow (\tl U,\tl U',\omega(0))$
 which is homotopic to a conjugacy
 $\psi : (U,U',\omega(0))\rightarrow (\tl U, \tl U', \omega(0))$
relative $(\di U,\di U',\omega(0))$. Note that $h$ conjugates 
$f: \omega(0)\cup \di U'\rightarrow \omega(0)\cup \di U$ and 
$\tl f: \omega(0)\cup \di \tl U'\rightarrow \omega(0)\cup \di \tl U$.
A qc map $h$ as above will be called a {\it Thurston conjugacy}.

\smallskip{\it Remark.} It is enough to assume that 
$h$ is homotopic to $\psi$ rel
postcritical sets. Then one can extend it to a qc map
 $U\rightarrow \tl U$ which 
is homotopic to $\psi$ rel the bigger set as required above.\smallskip
 
 The following
result comes from the work of Thurston (see \cite{DH2,McM1})
 and Sullivan (see \cite{MvS,S}). It originates the ``pull-back method" in 
holomorphic dynamics. 

\begin{lem}\label{pull-back}
If two quadratic-like maps are Thurston equivalent then they are qc equivalent.
\end{lem}

\begin{pf} We will use the notations for the domains and maps preceding the
statement of the lemma.  Let
$U^n$ be the preimages of $U$ under the iterates of $f$, and let $c=f(0)$. 
Let $h$ has dilatation $K$.
 
Since $h(c)=\tilde c$, we can lift $h$ to a $K$-qc map
$h_1: U^1=\tilde U^1$ homotopic to $\psi$ rel $(\di U^1, \di U^2, \omega(0))$.
(Note that the dilatation of $h_1$ is the same as the dilatation of $h$,
since the lift is analytic).  Hence
$h_1=h$ on these sets, and we can extend $h_1$ to
$U\sm U^1$ as $h$ (keeping the same notation $h_1$). By the Gluing Lemma
from the Appendix this extension has the same dilatation $K$. 
Moreover, this map is homotopic to
$\psi$ rel $(\omega(0), \cup_{1\leq k\leq 2} \di U^k)$. Also, it conjugates
 $f: \omega(0)\cup(U^1\sm U^2)\rightarrow\omega(0)\cup( U^0\sm U^1)$
to the corresponding tilde-map
(notice  that $h_1$ is a conjugacy on a bigger set than $h$).

Let us now replace $h$ with $h_1$ and repeat the procedure. We will construct a $K$-qc map
$h_2: U\rightarrow \tl U$ homotopic to $\psi$ rel 
$(\omega(0), \cup_{1\leq k\leq 3}\di U^k)$ and conjugating 
$f:\omega(0)\cup(U^1\backslash U^3)\rightarrow\omega(0)\cup( U\backslash U^2)$
to the corresponding tilde-map.

Proceeding in this
way we construct a sequence of $K$-qc maps $h_n$ homotopic to $\psi$ rel 
$(\omega(0), \cup_{1\leq k\leq n+1}\di U^k)$ and
conjugating  
$f: \omega(0)\cup(U^1\backslash U^{n+1})\rightarrow
\omega(0)\cup( U\backslash U^n)$ to the
corresponding tilde-map.
By the Compactness Lemma from the Appendix,
 we can select a converging subsequence
$h_{n(l)}\to h$. The 
limit map $h$ is a desired qc conjugacy.
\end{pf}

 The method  used in the above  proof  is called ``the pull-back argument". 
The idea is to
start with a qc map respecting some dynamical data, and then pull it back 
so that it will respect some new data on each step. In the end  it
becomes (with some luck)  a qc conjugacy. 
%This  idea will be heavily used in what follows.

{\it Remark.} For infinitely renormalizable maps of bounded type with 
{\it a priori} bounds, McMullen proved that the postcritical
set $\omega(0)$ has  bounded geometry \cite{McM3}. It easily follows that
there is a qc map $h: (\C, \omega(f,0))\ra (\C,\omega(\tl f,0))$
conjugating $f$ to $\tl f$ on their postcritical sets. This is close
to being a Thurston conjugacy but not quite the same, as $h$ may be in a 
wrong homotopy class.

\subsection{Approximating sequence of homeomorphisms}\label{strategy} 
So we need to construct a Thurston conjugacy. 
We will construct it as a limit of an appropriate sequence of
maps. Take a sufficiently small $\eps>0$, 
and consider the corresponding sequence
of standard neighborhoods $\U^m=\cup_i U^m_i\equiv U^m_i(\eps)$  
(see \secref{sec:standard neighborhoods}). 
By  \corref{exp decay of U} there is an $l$ such that $\U^m$ is well
 inside  $U^{m-l}$ (that is, the annulus $U^{m-l}\sm U^m$ has
a  definite modulus).
Moreover, by \lemref{omega removable} $\cap\J^m=\omega(0)$.

We will consecutively construct a sequence of homeomorphisms 
\begin{equation}\label{approx sequence}
h_m: (\C,\U^m,\J^m)\rightarrow (\C,\tl\U^m,\tl\J^m)
\end{equation}
such that

\smallskip\noindent (i) $h_0$ is a  topological conjugacy;

\smallskip\noindent (ii) $h_{m+1} $ is  homotopic to  $h_{m}$ rel 
$(\J^{m+1}\cup (\C\setminus \U^{m-l}))$. 
In particular  $h_{m+1}|\J^{m+1}=h_{m}|\J^{m+1}$ and
$h_{m+1}|(\C\setminus \U^{m-l})=h_m|(\C\setminus \U^{m-l})$.

\smallskip\noindent (iii) 
The $h_{m}$  are $K_*$-qc on $\U^{m-1}\sm \J^m$, with dilatation
 $K_*$ depending only on the choice of limbs and {\it a priori} bounds;

\smallskip\noindent (iv) $\Dil(h_{m}|U^{m-l})\leq 4K_*^4\Dil(h_{m-1}|U^{m-l}).$

Such a  sequence will do the job:

\begin{lem}\label{limit}
 A sequence $h_m$ satisfying the above three properties
  uniformly converges to a Thurston conjugacy.
\end{lem}

\begin{pf}
By the  second property, 
this sequence eventually stabilizes outside $\cap \J^m$ and
thus it pointwise converges to a homeomorphism
$h: (\C,\cap\J^m)\ra(\C,\cap\tl\J^m)$.
 By the last two properties, the dilatation of $h_m$  on  $\U^{m-l}\cap\J^m$
is at most $4^l K_*$. Hence $h$ is quasi-conformal on $\C\sm \cap \J^m$.
% By the Compactness Theorem, on $\C\setminus \cap\J^m$
%there is  a limit 
%$h=\lim h_{m(l)}: \C\setminus \cap\J^m\ra \C\setminus\cap\tl\J^m $,
% which is a qc homeomorphism. 
But by \lemref{omega removable} 
$\cap\J^m=\omega(0)$  is a removable Cantor sets.
 Hence $h$ admits a qc extension across $\omega(0)$. 

Further,  $h$ is homotopic to $h_0$ rel $\omega(0)$. 
Indeed, let $h^t$, $1-2^{-m}\leq t\leq 1-2^{-(m+1)}$,
 be a homotopy between $h_m$ and $h_{m+1}$
given by (iii). 
Let $\eps_m=\max_i\diam U^m_i$. As the $\U^m$ shrink to a  Cantor set,
$\eps_m\to 0$. As $h(U^{m-l}_i)= h^t(U^{m-l}_i)=\tl U^{m-l}_i$,
 $1-2^{-m}\leq t< 1$, 
 the uniform distance between $h$ and $h^t$ is at most $\eps_{m-l}$.
 It follows that 
the $h^t$ uniformly converge to $h$ as $t\to 1$. 
Hence $h$ is homotopic to $h_0$
rel $\omega(0)$. 

Since $h_0$ is a topological conjugacy by (i), $h$ is a Thurston conjugacy.
\end{pf}

\subsection{Construction of $h_0$}\label{sec: beginning}
Let us supply the exterior $\C\sm \cl\D$ of the unit disk,
with the hyperbolic metric $\rho$. 
The hyperbolic length of a curve $\gamma$ will be denoted by $l_\rho(\gamma)$,
while it Euclidean length will be denoted by $|\gamma|$. 

\begin{lem}\label{commuting maps} Let $A$ and $\tl A$
 be two (open) annuli whose inner boundary is the circle  $\T$.
 Let $\omega:A\ra \tl A$ be a homeomorphism commuting with 
$P_0: z\mapsto z^2$ near $\T$.
Then $\omega$ admits a continuous extension to a map
 $A\cup\T\ra \tl A\cup \T$ identical  on
the circle.
\end{lem}

\begin{pf} Given a set $X\subset A$, 
let $\tl X$ denote its image image by $\omega$.
 Let us take a configuration consisting of 
 a round annulus $L^0=\A[r,r^2]$ contained in
$A$, and an interval $I_0=[r, r^2]$.
Let $L^n=P_0^{-n}L^0$, and $I^n_{k}$ denote the components of $P_0^{-n} I^0$,
$k=0,1,\dots, 2^n-1$.
%Moreover, we count this components anticlockwise starting
%with $I_{n,0}\subset \R$. Note that
%$\cup_n I{n,0}\equiv \II_0$ is an  connected arc. 
The intervals $I^n_{k}$ subdivide
the annulus $L^n$ into $2^n$ ''Carleson boxes" $Q^n_{k}$.

 Since the (multi-valued) square root map
$P_0^{-1}$ is infinitesimally contracting in the hyperbolic metric,
 the hyperbolic diameters of the boxes
$\tl Q^n_{k}$ are uniformly bounded by a constant $C$.

Let us  now show that $\omega$ is a hyperbolic quasi-isometry near the circle,
 that is,
there exist $\eps>0$ and $A, B>0$ such that
\begin{equation}\label{quasi-isometry}
A^{-1}\rho(z,\zeta)-B\leq\rho(\tl z, \tl \zeta)\leq A\rho(z,\zeta)+B,
\end{equation}
provided $z,\zeta\in \A(1,1+\eps),\;  |z-\zeta|<\eps$.

Let $\gamma$ be the arc of the hyperbolic geodesic joining $z$ and $\zeta$.
 Clearly it is 
contained in the annulus $\A(1,r)$,
 provided $\eps$ is sufficiently small. 
Let $t>1$ be the radius of the circle 
$\T_t$ centered at 0 and tangent to $\gamma$. Let us
replace $\gamma$ with a {\it combinatorial geodesic} 
$\Gamma$ going radially up from $z$ to
the intersection with $\T_t$, 
then going along this circle, and then radially down
to $\zeta$. Let
$N$ be the number  of the Carleson boxes intersected  by $\Gamma$. 
Then one can easily see that
$$\rho(z,\zeta)=l_\rho(\gamma)\asymp \l_\rho(\Gamma)\asymp N,$$
provided $\rho(z,\zeta)\geq 10\log(1/r)$ (here $\log(1/r)$
 is the hyperbolic size of the boxes
$Q^n_k$).

On the other hand
$$\rho(\tl z, \tl \zeta)\leq l_\rho(\tl \Gamma)\leq CN,$$
so that $\rho(\tl z, \tl \zeta)\leq C\rho(z,\zeta)$, 
and (\ref{quasi-isometry}) follows.
 
But quasi-isometries of the hyperbolic plane admit
 continuous extensions to $\T$
(see, e.g., \cite{Th2}). 
Finally, it is an easy exercise to show that the only homeomorphism
of the circle commuting with $P_0$ if identical.
\end{pf}

\begin{lem}\label{commuting maps-1} Let $f$ be a quadratic-like map.
Let $A$ and $\tl A$ be two (open) annuli whose inner boundary is  $J(f)$.
 Let $\omega: A\ra \tl A$ be a homeomorphism commuting with $f$ near $J(f)$.
 Then $\omega$
admits a continuous extension to a map $A\cup J(f)\ra \tl A\cup J(\tl f)$ 
identical  on the Julia set.
\end{lem}

\begin{pf} 
By the Straightening Theorem, we can assume without loss of generality
 that $f=P_c: z\mapsto z^2+c$ is a quadratic polynomial. 
Let $R: \C\sm K(f)\ra  \C\sm \cl\D$
 be the Riemann mapping normalized by $R(z)\sim z$ near infinity. 
It conjugates
$P_c$ to $P_0: z\mapsto z^2$.

Let $\omega^\#=R\circ h\circ R^{-1}: \C\sm \cl D\ra \C\sm \cl D$. 
Then $\omega^\#$ commutes
with   with $P_0$ in an open annulus attached to the circle $\T$. 
By \lemref{commuting maps},
$\omega^\#$ continuously extends to $\T$ as $\id$. Hence for  any
$\eps>0$ there is an $r>1$  such that $|\omega^\#(z)-z|<\eps$
for $z\in \A(1,r)$. 

Let us show that the hyperbolic distance $\rho(\omega^\#(z), z)$
 is bounded if $|z|<2$. Clearly 
$\rho(\omega^\#(z), z)\leq  C(r)$, provided $1<r\leq |z|<2$. 
Let $r^{\frac 12} \leq |z| \leq r$, $\zeta=\omega^\#(z)$.
%It is easy to see that for sufficiently
%small $\eps>0$, 
%there is a curve  $\gamma\subset E_{r+\eps}$ 
%joining $P_0(z)$ and $P_0(\zeta)$, with  
%$|\gamma|< 2\eps$ 
%(one can take an arc of the circle centered at 0 and a straight interval).
Let us consider the hyperbolic geodesic $\gamma$ joining $z$ and $\zeta$.
Clearly $|\gamma|<O(\eps)$.
 Then
$P_0^{-1}\gamma$ consists of two  symmetric curves $
\sigma$ and $-\sigma$ of Euclidean  length
$O(\eps)$. 
One of these curves, say $\sigma$, 
joins $z$ with a preimage $u$ of $P_0(\zeta)$.  Then 
$|z+u|>2-O(\eps)>\eps$, so that $-u\not=\zeta$. Thus $u=\zeta$.

As the square root map $P^{-1}_0$ is infinitesimally 
contracting in the hyperbolic metric,
$$\rho(z,\zeta)\leq l_\rho(\sigma)\leq 
l_\rho(\gamma)=\rho(P_0(z), P_0(\zeta))\leq C(r).$$

Take now any point $z$  in the annulus 
$\A(r^{\frac 1 4}, r^{\frac 1 2}\}$. Using the same argument
we conclude that $\rho(z, \omega^\#(z))\leq C(r)$ (with the same $C(r)$. 
 By induction,  the same bound holds for all $z$.

Now we can complete the proof. 
Since the Riemann mapping $R$ is a hyperbolic isometry,
the hyperbolic distance between $\omega(z)$ 
and $z$ in $\C\sm J(P_c)$ is also bounded
near $J(P_c)$. Hence the Euclidean distance 
$|z-\omega(z)|$ goes to 0 as $z\to J(f)$.
It follows that the extension of $\omega$ as the identity on the Julia set
is continuous.
  \end{pf}

\begin{cor}\label{matching} Let $f$ and $\tl f$ be two topologically 
equivalent quadratic-like map, 
and let $\psi$ be a topological conjugacy between them.
Let $A$ and $\tl A$ be two open annuli whose inner boundaries are  $J(f)$
and $J(\tl f)$ respectively.
 Let $h: A\ra \tl A$ be a homeomorphism conjugating $f$ and $\tl f$ on these
annuli. Then
$h$  matches with $\psi$ on the Julia set, that is $h$
admits a continuous extension to a map
$A\cup J(f)\ra \tl A\cup J(\tl f)$  coinciding with $\psi$  on the Julia set.
\end{cor}

\begin{pf} Apply \lemref{commuting maps-1} to the homeomorphism 
$\omega=\psi^{-1}\circ h$
commuting with $f$.
\end{pf}

\begin{lem}[\cite{DH2}] \label{beginning}
 If quadratic-like maps  $f$ and $\tl f$ are
topologically  conjugate then there is conjugacy 
$h_0$ which is quasi-conformal outside the
Julia sets.
\end{lem} 

\begin{pf}  Given an annulus $A$, let
$\di_o A$ and $\di_i A$ stand for its outer and inner boundary components. 
Let us select a closed fundamental annulus $A$ for $f$ with smooth boundary, 
and let 
$A^n=f^{-n}A$. Let $\tl A$ and $\tl A^n$ be similar objects for $\tl f$. 
Then there is a
diffeomorphism $\phi: A\ra \tl A$ such that 
$$\phi(fz)=\tl f(\phi z), \quad z\in \di_i A.$$

This diffeomorphism can be lifted to a diffeomorphism 
$\phi_1: A^1\ra \tl A^1$ with the same
qc dilatation and such that
$$\phi_1(z)=\phi(z),\; z\in \di_o A^1,
\quad and \quad \phi_1(fz)=\tl f(\phi_1 z),
 \; z\in \di_i A^1.$$
In turn, $A^1$ can be lifted to a diffeomorphism 
$\phi_2: A^2\ra A^2$ with the same
dilatation, which matches with $A^1$ on $\di_o A^2$ 
and respects dynamics on $\di_i A^2$, etc.
    
By the Gluing Lemma from the Appendix, 
these diffeomorphisms glue together into a single
quasi-conformal map $h_0: A\sm J(f)\ra \tl A\sm J(\tl f)$
 conjugating $f$ and $\tl f$.
 
On the other hand, 
let $\psi$ be a topological conjugacy between $f$ and $\tl f$ near the
Julia sets. Then by \corref{matching},  $h_0$ matches with $\psi$ on
$J(f)$.
\end{pf}
 
\subsection{Adjustment of $h_m$}\label{adjustment sec}
\comm {Let us first fix some terminology and notations.
Recall that $\Dil(h)$ stands for the dilation of a qc map $h$. 
Given a conformal annulus $A$, let $R: A\ra A^u$ be the Riemann mapping onto
the standard cylinder $A^u=\T\times [0, H]$ where $\T$ is the unit circle.
It endows $A$ with the {\it cylinder coordinates} (angle, height). The {\it horizontal}
and {\it vertical} curves on $A$ are the cylinder coordinate lines, that is,
the $R$-preimages of the horizontal circles and vertical intervals in $A^u$.
A map $\phi: A\ra\tl A$ between two conformal annuli is called {\it affine}
if it becomes affine after uniformization.}

Recall that $p_m$ is the period of the little Julia sets $J^m_j$, and 
$F_m=f^{p_m}$ is the  corresponding quadratic-like map near $J^m_j$.  
Let  $\U^m=\cup U_j^m$ be
a standard neighborhood of the little Julia orbit $\J^m=\cup B^m_j=\cup J^m_i$, 
with a definite space in between the $ U^m_j$ and definite
annuli $U^m_j\setminus B^m_j$, and let
 $S_m: \U^m_j\sm \B^m_j \ra \A(1,4) $ be the standard
straightenings  (\ref{standard straightening}). Its dilatation is 
bounded by a constant  $K_*$ depending only on the choice of secondary limbs 
and {\it a priori} bounds.
%The  annuli $U^m_j\sm B^m_j$ have the standard coordinates
%coming from the straightening of $F_{m-1}$ with a definite dilation. 
Let $U^m_j(t)=S_m^{-1}\A(1,t)$  (note that $U^m_j\equiv U^m_j(4)$).
 The notation $\U^m(t)$ is self-evident.

 We say that a homeomorphism
 $\phi: U^m_j(2)\sm B^m_j\ra \tl U^m_j(2)\sm \tl B^m_j$ 
is {\it standard} near the bouquet $B^m_j$ if it  is identical 
 in the standard coordinates on  $U^m_j(2)$, that is,
\begin{equation}\label{standard map}
\tl S_m\circ \phi|U^m_j(2)=S_m.
\end{equation}
The dilatation of  such a map is bounded
by $K_*^2$.  
Note also that by \corref{matching}, the standard map 
 admits a homeomorphic extension  across the Julia bouquet.

% Let $\xi>0$ be a lower bound for the moduli of 
%the above collars and the space,
%while $1/\xi$ be an upper bound for the $\mod (T^m_j\sm B^m_j)$.   

%Let $\U^m\subset \U^{m-l}$ (see \lemref{inclusion}). 
We will now adjust the map
$h_m$ so that it will become standard near $\J^m$. 

\begin{lem}\label{adjustment} Take an $l$ as in \secref{strategy}.
Let a homeomorphism $h_m: (\C,\J^m)\ra (\C,\tl\J^m)$
 be a conjugacy on $\J^m$ 
 and be $K_m$-qc on $\U^{m-l}\sm \J^m$. 
 Then there is a homeomorphism 
$$\hat h_m: (\C,\U^m,\J^m)\ra (\C, \tl\U^m \tl\J^m)$$ 
 homotopic to $h_m$ rel $(\J^m\cup (\C\sm \U^{m-l}))$, such that
$\Dil(\hat h^m|(\U^{m-1}\setminus \J^m))\leq 4K_*^4\cdot K_m,$
 and $h_m: \U^m(2)\sm \J^m\ra \tl \U^m(2)\sm \tl \J^m $ is standard.
% $-qc conjugacy between $f$ and $\tl f$. 
\end{lem}

\begin{pf}
Let us consider a retraction 
$\psi^t_j: U_j(4)\setminus B_j\ra U_j(4)\sm U_j(2)$ 
which is the affine vertical  contraction in the standard coordinates.  
Its dilation is bounded by $2K_*^2$.   
 Let us extend the $\psi^t_j$ to a homeomorphism
 $\psi:\C\sm \J\ra \C\sm U(t)$ by identity on $\C\sm \U(4)$.
 By the Gluing Lemma from the Appendix, 
$\psi$ is also  $2K_*^2$-qc. 

Let us now define a homeomorphism 
$ h^t: (\C, \U^t, \J)\ra (\C, \tl\U^t, \tl\J)$ as follows:
 $$ h^t|(\C\sm \U^t)=\tl\psi^t\circ h\circ(\psi^t)^{-1},$$
while $h^t: \U^t\ra\tl \U^t$ is standard. Then $h^1$ is a desired adjusted map
(homotopic to  $h^0=h$ via the $\{h^t\})$.  
\end{pf}

In what follows we will assume that $h_m$ is adjusted as in \lemref{adjustment},
and will skip the ''hat" in the notation for the adjusted map. 

\subsection{Beginning of the construction of $h_{m+1}$}\label{beginning sec}
% Let us go back to  the situation described by 
%formula (\ref{inclusions}).
 Let $p_m$ denote the
combinatorial rotation number of the  $\alpha$-fixed of the Julia sets $J^m_i$.
 Consider
the configurations $\RR^m_i$ of $2p_m$ rays   landing at the 
 $\alpha$-fixed and co-fixed points of the $J^m_i$.
 Let $\Omega^{m,0}_s\equiv \Omega^m_s$ stand for
the component of $U^m_j\sm \RR^m_i$ containing $J^{m+1}_s$, and let 
$\Omega^{m,1}_s\subset \Omega^m_s$ 
be the component of $F_m^{-p_m}\Omega^m_s$ containing $J^{m+1}_s$,
so that 
\begin{equation}\label{double covering }
G_m\equiv F_m^{\circ p_m}: \Omega^{m,1}_s \ra \Omega^m_s
\end{equation} 
is a double branched covering.
 The boundaries of these domains 
are naturally marked with the standard coordinates.
({\it Marking} of a curve  means its preferred parametrization.)
 As the map 
$$ h_m: (\C,  U^m_j,  \Omega^m_s, \Omega^{m,1}_s, J^{m+1}_s)\ra
    (\C, \tl U^m_j, \tl\Omega^m_s, \tl\Omega^{m,1}_s, \tl J^{m+1}_s)$$
is standard on  the $U^m_j$, it respects this marking.

Since the configurations $(\cup \RR^m_s, \cup \di \Omega^{m,1}_s)$
 have bounded geometry
(see \S 4 of Part I), there is a 
qc map with a bounded dilatation
\begin{equation}\label{initial map}
\Psi_{m}: (\C, U^m_j, \Omega^m_s,\Omega^{m,1}_s)\ra 
(\C, \tl U^m_j, \tl\Omega^m_s,\tl\Omega^{m,1}_s)
\end{equation}
coinciding with $h_m$ on $\C\sm \Omega^m$ and
respecting the boundary marking  (in particular, it conjugates
$F_m: \di\Omega^{m,1}_s\ra\di\Omega^m_s$ and 
$\tl F_m: \di\tl\Omega^{m,1}_s\ra\di\tl\Omega^m_s$).
Moreover $\Psi_m$ is homotopic to $h_m$ rel 
$((\C\sm\U^m)\cup\di\Omega^m_s\cup\di\Omega^{m,1}_s)$,  since all
 regions complementary to this set are simply connected Jordan domains.

Note however that unlike $h_m$,
the map $\Psi_m$ does not respect dynamics on the little Julia
sets. We need to pay temporarily this price in order to make $\Psi_m$
 globally quasi-conformal.

\subsection{Construction of $h_{m+1}$ in  the immediately renormalizable case}
\label{immediate case}
 Let us consider the double covering  (\ref{double covering }).
In the immediately renormalizable case, 
$$G_m^{\circ n}0\in \Omega^{m,1}_s,\; n=0,1,2,\dots$$
Moreover, there is a nest of topological disks
$$\Omega^{m,0}_s\supset \Omega^{m,1}_s\supset \Omega^{m,2}_s\supset\dots$$
shrinking to the little Julia set $J^{m+1}_s$, and such that
$G_m: \Omega^{m,n}_s\ra \Omega^{m,n-1}_s$ is a branched double covering.
 The complement $Q^{m,n}_s=\Omega^{m,n-1}_s\sm \Omega^{m,n}_s$
 consists of $2^n$ quadrilaterals.  

As $G_{m}: Q^{m,n}_s\ra Q^{m,n-1}_s$ is an unbranched covering,
the map $\Psi: Q^{m,1}_s\ra \tl Q^{m,1}_s$ can be lifted to a qc map
$$\Psi_{m,n}: Q^{m,n}_s\ra \tl Q^{m,n}_s$$ with the same dilatation
homotopic to $h_m$ rel the boundary. 
Hence all these maps glue together in a single qc map with the same dilatation 
\begin{equation}\label{(m+1)t map}
h_{m+1}: \Omega^{m}_s\sm J^{m+1}_s\ra \tl \Omega^{m}_s\sm \tl J^{m+1}_s
\end{equation}
 equivariantly homotopic to $h_m$ rel $\cup_n\di \Omega^{m,n}$. 

Let $\psi^t$ be the corresponding homotopy, and $\rho$ be the
hyperbolic metric in $\tl \Omega^{m}\sm\tl J^{m+1}$. Then by equivariancy
$\rho(\psi^t(z), h_{m}(z))\leq C$. 
Hence $|\psi^t(z)-h_{m}(z)|\to 0$ as $z\to J^{m+1}_s$
uniformly in $t$. It follows that the homotopy $\psi^t$ 
 extends    across the little Julia set $J^{m+1}_s$.
Thus the map  (\ref{(m+1)t map}) is extends across 
$J^{m+1}_s$ to a homeomorphism homotopic
to $h_m$ rel ($\di \Omega^m_s\cup J^{m+1}_s$).     

Outside the  $\bigcup_s\Omega_s^m$ let $h_{m+1}$ coincide with $h_m$.
This provides us with the desired map $h_{m+1}$.

\section{Through the principal nest}\label{sec: principal nest} 

{\it In what follows we will assume that $R^m f\equiv F_m$ 
is not immediately renormalizable.}

\subsection{Teichm\"uller distance between the configurations of puzzle pieces}
\label{sec: main lemma}
Let us make a choice of a standard neighborhood $U^m$ of the Julia bouquet $B^m$
and the corresponding  standard straightening  $S_m$,
see (\ref{standard straightening}). When $F_{m-1}$ is not immediately
renormalizable,
this provides us with a family $\YY$  of puzzle pieces $Y^{(k)}_i$,
see \S 2.6 of Part I. 

In the immediately renormalizable case let us start 
the puzzle in a slightly different way.
 Namely, let us consider a {\it degenerate}
domain of $F_m$ (see \S 2.5 of Part I)
%instead of taking a neighborhood
%of the Julia set $J^m$ bounded by an equipotential,
  bounded by external rays landing at  fixed and co-fixed
points $\alpha_{m-1}\equiv \beta_m$, $-\beta_{m}$, and two pieces of
standard  equipotentials of $F_{m-1}$.  Then play the puzzle by cutting
this domain with external rays of $F_{m-1}$ landing at $\alpha_m$, and pulling
them back. One can easily see that this beginning is equally suitable
for the puzzle game as the usual one.
%(the only  place which should be actually checked is that the first piece $V^0$
%of the principal nest is  contained in $\inter Y^{(0)}$, see \S 3.2 of Part I).

As  the puzzle pieces $Y^{(k)}_i$ are bounded by
equipotentials and rays, they  bear the {\it standard boundary marking},
i.e. the parametrization $S_m^{-1}$ by the corresponding straight intervals
or circle arcs. 

Since $h_m: U^m\ra \tl U^{m}$ is the standard conjugacy
(see (\ref{standard map})), 
it maps  the pieces $Y^{(k)}_i$ to the corresponding tilde-pieces 
 $\tl Y^{(k)}_i$ respecting the boundary marking.
Given   some family of puzzle pieces 
$P_i\in \YY$ contained in some $Y\in \YY$, let us say that a homeomorphism
$$\phi: (Y, \cup P_i)\ra (\tl Y, \cup \tl P_i)$$ 
is a {\it pseudo-conjugacy} if it is
homotopic to $h_m$ rel the boundary
$(\di Y,\cup\di P_i)$.
Note that if $f^l : P_i\ra Y$ (or $f^l : P_i\ra P_j$) for  some iterate
of $f$ and
some puzzle pieces of our family, then the pseudo-conjugacy $\phi$
is a true conjugacy between the boundary maps
$f^l: \di P_i\ra \di Y$ and  $\tl f^l|\di\tl P_i\ra \di\tl Y$
(correspondingly $\di P_j$ instead of $\di Y$).

In particular,  the above terminology will be applied to
 the principal nest of puzzle pieces  (see  \S 3 of Part I):
\begin{equation}\label{principal nest}
Y^{(m,0)}\supset V^{m,0}\supset V^{m,1}\supset\dots,\qquad 
V^{m,n}_0\equiv V^{m,n}, \qquad \cap_n V^{m,n}=J^{m+1},
\end{equation}
and the corresponding generalized renormalizations
   $g_{m,n}: \cup_i V^{m,n}_i\ra V^{m,n-1}.$ 
%$V^{m,n}_0\equiv V^{m,n}$.  
%The puzzle pieces $V^{m,n}_i$ are bounded by
%equipotentials and rays of the quadratic-like map $F_m$ 
%and hence bear the marking
%coming from the straightening of this map on the standard neighborhood $U^m$.
%As $h_m: U^m\ra \tl U^{m}$ is the standard conjugacy, 
%(i.e., it comes from the straightenings), 
%it maps $(V^{m,n+1},\cup_i V^{m,n}_i)$ to 
% $(\tl V^{m,n+1},\cup_i\tl V^{m,n}_i)$ respecting the boundary marking.
%So we say that 
%\begin{equation}\label{pseudo-conjugacy def}
%\phi:(V^{m,n+1},\bigcup_i V^{m,n}_i)\ra 
%(\tl V^{m,n+1},\bigcup_i\tl V^{m,n}_i)
%\end{equation}
%is a {\it pseudo-conjugacy} (between $g_m$ and $\tl g_m$) if it is
%homotopic to $h_m$ rel the boundary
%$(\di V^{m,n+1},\cup_i\di V^{m,n}_i)$. 
%In particular, it is a true conjugacy
%between the boundary maps.

{\it Teichm\"ullere distance}
 $\dist_T$ between $(V^{m,n+1},V^{m,n}_i)$ and 
 $(\tl V^{m,n+1},\tl V^{m,n}_i)$ is defined as
 $\inf_\phi\log K_\phi$ as $\phi$ runs over all qc pseudo-conjugacies
$(V^{m,n+1},\cup_i V^{m,n}_i)\ra (\tl V^{m,n+1},\cup_i\tl V^{m,n}_i)$.  

\proclaim Main Lemma [\cite{L5}, \S 4)]. 
The configurations $(V^{m,n+1},V^{m,n}_i)$ 
and  $(\tl V^{m,n+1},\tl V^{m,n}_i)$
stay bounded Teichm\"uller distance away (independently of $m$ and $n$).

The rest of this section, except the final subsection,\secref{spreading}), 
will be occupied with the proof of this lemma which follows \cite{L5}, \S 4.
As the level $m$ is fixed, {\it in what follows we will  skip the label $m$
in the notations of $V^{m,n}_i\equiv V^n_i$, $g_{m,n}\equiv g_n$ etc. } 
(unless it may  lead to a confusion). 
In what follows referring to a qc-map,
 we will mean that it has a definite dilatation
(depending only on the selected limbs and {\it a priori} bounds).
% When we say that some
%quantity (modulus, dilatation etc.) is bounded,
% we also mean that the bound depends only
%on the selected limbs and {\it a priori} bounds.

\subsection{A point set topology lemma} 
In the statement below,
 the objects involved need not have any dynamical meaning. 

\begin{lem}\label{homotopy}
 Let $P_i$ be a family of  closed Jordan disks with disjoint
interiors contained in a domain $Y$,
 such that $\diam P_i\to 0$. Let $\tl P_i$, $\tl Y$ 
be  another family of disks with the
same properties. 

\noindent$\bullet$  Let $h: (Y,\cup P_i)\ra (\tl Y, \cup \tl P_i)$ 
be a one-to-one map,
which is a homeomorphism on $\cup P_i$ and on $X\equiv Y\sm (\cup \inter P_i)$.
 Then $h $ is a homeomorphism.
 
\noindent$\bullet$ 
 Let $h^i: (Y, \cup P_i)\ra (\tl Y, \cup\tl P_i) $, $i=0,1$, 
be two homeomorphisms
coinciding on $Y\sm \cup\inter P_i$.  
Then $h^i$ are homotopic rel $Y\sm \cup\inter P_i$.
\end{lem}

\begin{pf} Given an $\eps>0$, there exists an $N$ such that
$\diam(\tl P_n)<\eps$ for $n> N$.   Let $T=\cup_{1\leq i\leq N} P_i$. 
 Note that $h$ is 
continuous on
$X\cup T$. 

Given a point  $z\in Y$, let us show  that $h$ is continuous at it.
 This is certainly true
if $z\in \cup \inter P_i$, so let $z\in X$. We will show that
\begin{equation}\label{continuity} |h(z)-h(\zeta)|< 2\eps
\end{equation} for any nearby point $\zeta\in Y$. 
 Indeed, if $\zeta\in X\cup T$ it follows
from the above remark. Otherwise $\zeta\in P_j$ for some $j> N$,
 and there is  point  $u\in [z,\zeta]\cap\di P_j$. Then
$$|h(z)-h(\zeta)|\leq |h(z)-h(u)|+|h(u)-h(\zeta)|.$$  If $\zeta$ is sufficiently close to
$z$ then  the first term is at most $\eps$ by continuity of $h|X$. As the  second term  is
bounded by $\diam(P_j)<\eps$, and (\ref{continuity}) follows. 

Let us now prove the second statement. As each
$P_i$ is simply connected, $h^0|P_i$  is homotopic to $h^1|P_i$ rel $\di P_i$. Let 
$h^t:\cup P_i\ra \tl P_i$ be a corresponding homotopy.  Extend it to the whole domain
$Y$ as $h^0$. We should check that this extension
 $h^t(z): (Y,\cup P_i)\ra (\tl Y, \cup \tl P_i)$ is continuous in two variables.  

Note first that for  $z\not\in\cup_{1\leq i\leq N} P_i\equiv T$,
\begin{equation}\label{est0} |h^t(z)-h^0(z)|<\eps.
\end{equation} 
\comm{As $h$ is a homeomorphism and $h^t: T\ra T$ is a homotopy, there is a $\delta>0$
such that 
$|h(z)-h(\zeta)|<\eps/3$ if $|z-\zeta|<\delta$, and
 \begin{equation}\label{est1} 
|h^t(z)-h^\tau(z)|,\eps/3 \quad if quad  z,\zeta\in T,\;
|z-\zeta|<\delta,\;|t-\tau|<\delta.
\end{equation}}
 Given a pair $(z,t)$, we will show that $|h^t(z)-h^\tau(\zeta)|<3\eps$ as
$(\zeta,\tau)$ is sufficiently close to $(z,t)$. 
To this end let us consider a few cases:

\noindent$\bullet$ If $z\in\inter \cup P_i$,
 it is true since $h^t|P_i$ is a homotopy.

\noindent$\bullet$ If $z,\zeta\in T$, it is true since $h^t|T$ is a homotopy.

\noindent$\bullet$ If $z\in\di T$ but $\zeta\not\in T$,
 then for $\zeta$ sufficiently close to $z$,
$$|h^t(z)-h^\tau(\zeta)|=|h^0(z)-h^\tau(\zeta)|\leq
|h^0(z)-h(\zeta)|+|h^\tau(\zeta)-h^0(\zeta)|<2\eps$$ 
by continuity of $h$ and (\ref{est0}).

\noindent$\bullet$ Let $z\not\in T$.
 Then sufficiently close points $\zeta$ don't belong to $T$ either.
 Hence by (\ref{est0})
and continuity of $h$,
$$|h^t(z)-h^t(\zeta)|\leq
|h^0(z)-h^0(\zeta)|+|h^t(z)-h^0(z)|+|h^t(\zeta)-h^0(\zeta)|<3\eps.$$ 
\end{pf}

\subsection{Expanding sets}
 Let us consider Yoccoz puzzle pieces $Y^{(N)}_i$ of depth
$N$  (see \S 2.6 of Part I), and let
$\YY^{(N)}$ denote the family of puzzle pieces $Y^{(N+l)}_j$ such that
$$f^k  Y^{(N+l)}_j\cap Y^{(N)}_0=\emptyset, \; k=0\dots,l-1 .$$
Let $K^{(N)}=\{z: F^k z\not\in Y^{(N)},\; k=0,1\dots\}$
Recall that an invariant set $K$ is called {\it expanding} 
if there exist constants
$C>0$ and $\rho\in (0,1)$ such that
$$|D F^k(z)|\geq C\rho^k,\; z\in K,\; k=0,1,\dots$$

\begin{lem}\label{shrinking} For a given $N$,
 $\diam Y^{(N+l)}_s\to 0$ 
as $Y^{(N+l)}_s\in \YY^{(N)}$ and $l\to\infty$. Moreover, 
the set $K^{(N)}$ is expanding.
\end{lem}

\begin{pf} Let us consider thickened puzzle pieces $\hat Y^{(N)}_i$ as in
Milnor \cite{M1} or \S 2.5 of Part I. 
Then $\inter(F \hat Y^{(N)}_i)$ contains $\hat Y^{(N)}_j$
 whenever $F Y^{(N)}_i\supset Y^{(N)}_j $  (recall
that the $Y^{(N)}$ are closed). Hence the inverse map 
$F^{-1}: \hat Y^{(N)}_j\ra \hat Y^{(N)}_i$
is  contracting by a factor $\lambda<1$ 
in the hyperbolic  metrics of the pieces under
consideration.

Let $Y^{(N+l)}_s\subset Y^{(N)}_i$. It follows
that the hyperbolic diameter of $Y^{(N+l)}_s$ in $Y^{(N)}_i$ is at most
$\lambda^{l}$, and the statement follows.
\end{pf}

\subsection{First landing maps}\label{sec: first pull-back}

Let us have a family of  puzzle pieces $P_i$ with disjoint interiors
contained in a 
puzzle piece $X$, where as usual $P_0\ni 0$ stands for the critical 
puzzle piece. Let us also have  a Markov map $G: \cup P_i\ra X$
which is univalent on all non-critical pieces  $P_i$, $i\not=0$,
 and the double branched
covering on the critical one, $P_0$. The Markov 
property means that if $\inter(G P_i\cap P_j)\not=\emptyset$ then
$G P_i\supset P_j$. Let $A$ be the corresponding Markov matrix:
$A_{ij}=1$ if $\inter(G P_i\cap P_j)\not=\emptyset$, and 
$A_{ij}=0$ otherwise.

Let $P\equiv P^0$. A string of labels
$\bar i=(i(0),\dots, i(l-1))$   is called {\it admissible} if
 $A_{i(k), i(k+1)}=1$ for
$k=0, \dots, l-2$), and  $i(k)\not=0$ for $k<l-1$. 
The length $l$ of the string will be denoted by  $|\bar i|$.
To any admissible string corresponds a
 {\it cylinder} of rank $l$ defined by the following property:
\begin{equation}\label{cylinders}
G^k P_{\bar i}^l\subset P_{i(k)},\; k=0,\dots l-2,\;\;
 G^{l-1} P_{\bar i}^l=P_{i(l-1)}.
\end{equation}
Note that   $G^{l-1}$ univalently maps $P_{\bar i}^l$ onto $P_{i(l-1)}$.

 Let us denote
by $\Omega_{\bar i}\equiv P_{\bar i}^l$
the cylinders mapped onto the critical puzzle piece (so that $i(l-1)=0$).
The {\it first  landing map} 
\begin{equation}\label{first landing}
T: \cup \Omega_i\ra P_0
\end{equation}
is 
defined as follows: $Tz=G^{l-1}z$ for $z\in \Omega_{\bar i}$, $|\bar i|=l$.
This map is univalent on all pieces $\Omega_i$ (identical on the
critical piece $\Omega_0$).

\begin{lem}\label{simple pull-back}
Let us  have a  $K$-qc pseudo-conjugacy
$H: (X, \cup P_i)\ra (\tl X, \cup \tl P_i)$ between $G$ and $\tl G$. 
Then there is  a $K$-qc pseudo-conjugacy 
$\phi: (X, \cup \Omega_j)\ra (\tl X, \cup \tl \Omega_j)$
which conjugates  the first landing maps $T$ and $\tl T$.
\end{lem}

\begin{pf}  Let us pull $H$ back to the pieces $P_i$,  $i\not=0$,
that is, let us consider the map 
$$H_1:  (P_i, \bigcup_j P^l_{ij})\ra (\tl P_i, \bigcup_j \tl P^l_{ij})$$ 
such that $\tl G\circ H_1|P_i=h\circ G|P_i$. 
Since $H$ is a pseudo-conjugacy, $H_1$ matches with $H$ on
$\cup_{i\not=0}\di P_i$. Hence these maps glue together into a single
map $K$-qc map equal to $H_1$ on $\cup\di P_i$, and equal to $H$ outside
of it. We will keep notation $H_1$ for this map.

Let us do the same pull-back with $H_1$. We will obtain a 
$K$-qc pseudo-conjugacy
$$H_2: (P, \cup P_i^1, \cup P_{ij}^2, \cup P_{ijk}^3)\ra
(\tl P, \cup \tl P_i^1, \cup \tl P_{ij}^2, \cup \tl P_{ijk}^3).$$
Repeating this procedure over again,
 we obtain a sequence of $K$-qc pseudo-conjugacies
$$H_s: \bigcup_{l\leq s}\bigcup_{ |\bar i|=l} P^l_{\bar i}\ra 
 \bigcup_{l\leq s}\bigcup_{ |\bar i|=l} \tl P^l_{\bar i}.$$
By the Compactness Lemma from the Appendix we can pass to a limit
$K$-qc map
 $$\phi: \bigcup_{l, \bar i} P^l_{\bar i}\ra 
 \bigcup_{l, \bar i} \tl P^l_{\bar i}.$$
By \lemref{homotopy} this map is homotopic to $h$ rel $(\di X\cup\di \Omega_j)$,
 and hence is a desired pseudo-conjugacy.
\end{pf}

Let us now do a bit more (assuming a bit more).
Let us consider the generalized renormalization of $G$ on $P_0$, that is, the 
first return map $g: \cup V_j\ra P_0$. Let $b= g(0)=G^t0$ be its critical value.

\begin{lem}\label{simple pull-back-2}
Let us have two $K$-qc pseudo-conjugacies
 $H_0: (X,\cup P_i)\ra (\tl X, \cup \tl P_i)$
and  $H_1: (P_0,b)\ra (\tl P_0,\tl b)$.
Then there exist a $K$-qc pseudo-conjugacy 
$\psi: (P_0, \cup V_i)\ra (\tl P_0, \cup \tl V_i)$ between $g$ and $\tl g$.
\end{lem}

\begin{pf} As $H$ and $H'$ match on $\di P_0$, 
they  glue together into a singe $K$-qc pseudo-conjugacy 
$H: (X, \cup P_i, b)\ra (\tl X, \cup\tl P_i, \tl b)$
coinciding with $H_1$ on $P_0$ and coinciding with $H_0$ on $X\sm P_0$
(see the Gluing Lemma in the Appendix).
By \lemref{simple pull-back}, there is a
$K$-qc   map
$\phi:  (X, \cup \Omega_j)\ra (\tl X, \cup \tl \Omega_j)$
homotopic to $h$ rel ($\di X\cup\di\Omega_j$), and 
conjugating the first landing maps. As $H: b\mapsto \tl b$, we have: 
$\phi: G^k 0\mapsto \tl G^k 0$, $k=1,\dots, t$. In particular, $\phi$ respects the
$G$-critical values: $G(0)\mapsto\tl G(0)$.

Recall that the domains $V_i$ are the pull-backs of the $\Omega_j$ by
$G: P_0\ra X$, that is, the components of $(G|P_0)^{-1}\Omega_j$. 
It follows that $\phi$ can be lifted to a $K$-qc map
$\psi: (P_0, \cup V_i)\ra (\tl P_0, \cup\tl V_i)$ homotopic to $h$ rel 
$(\di P_0\cup\di V_i)$.
 (This  lift is uniquely determined by the diagram $\tl G\circ\psi|P_0= \phi\circ G|P_0$
and the homotopy condition.)

 This  map $\psi$ is the desired pseudo-conjugacy.
\end{pf}

\subsection{Initial constructions}\label{initial construction}
Now the reader should consult \S 3.2 of Part I of this paper \cite{L6}, 
where the initial Markov partition  (3-3) of the Yoccoz  puzzle piece
 $Y^{(0)}$ is constructed.
We will apply it to the renormalized map $F$. Let us recall some notations.
The first piece of the partition, $Y\equiv Y^{(0)}$,
 is bounded by the external rays landing
at the fixed point $\alpha$, and the equipotential $E$.
The central piece of this partition,
 $V^0$, is the first piece of the principal nest.
It is obtained by pulling back a puzzle piece $Z_\nu^{(1)}$
 attached to the co-fixed point
$\alpha'$ (that is, $F(\alpha')=F(\alpha)$).
 There is a double branched covering $F^{s}: V^0\ra Z_\nu^{(1)}$. 
All the puzzle pieces  of the initial partition 
intersecting the Julia set $J(F)$ are univalent pull-backs of either
$Y$ or $V^0$. Let us denote the pieces of this partition by $P_i$, 
in such a way that
$P_0\equiv V^0$, $P_i\equiv Z^{(1)}_i$, $i=1,\dots p-1$, where $p$ 
is the number of
external rays of $F$ landing and $\alpha$. With these notations,
\begin{equation}\label{initial partition}
Y\cap J(F)=\bigcup (P_i\cap J(F))\cup K,
\end{equation}
where $K$ is the  residual Cantor set 
(of the points whose orbits never land at
$\cup_{0\leq i\leq p-1} P_i$). 

\begin{lem}\label{initial removability}
In the decomposition (\ref{initial partition}), 
$\diam P_i\to 0$ and the set $K$ is a
removable Cantor set.
\end{lem}

\begin{pf} \comm{Let us partition a neighborhood of the Julia set by the 
puzzle pieces $W_j$ of the same {\it depth}  $N$ 
as $V^0$  (see \S 2.6 of Part I), where $W_0\equiv V^0$. 
Let us consider thickened pieces $\hat W_i$ 
as in  Milnor \cite{M1} or \S 2.5 of Part I. Then $\inter(F^p \hat W_i)$
contains $\hat W_j$ whenever $F^p W_i\supset W_j$ 
(recall that $W_i$ are closed).
Hence the inverse map $F^{-p}: \hat W_j\ra \hat W_i$ is  contracting
by a factor $\lambda<1$ in the hyperbolic  metrics of the pieces under
consideration.
       Given a puzzle piece $P_i$ of
depth  $N+l+1$ we have:
$$F^{pk} P_i\subset W_{k(i)} \;\; with \; k(i)\not=0, \;\; k=0,\dots, l .$$
It follows that the hyperbolic diameter of $P_i$ in $\hat W_{i(0)}$ is at
most $\lambda^{l}$, and the first statement follows.   end comment}
The first statement follows from \lemref{shrinking}.
To prove removability of $K$, let us consider the domains $Q_1$ and $Q_2$
introduced in \S 3.2 of Part I. Then $F^p Q_i$ covers $Q_1\cup Q_2$,
and $K$ is the  set of points which never escape $Q_1\cup Q_2.$ By a 
little thickening of these domains, we obtain a Bernoulli map
$F^p: \hat Q_1\cup \hat Q_2\ra \C$ 
(so that $\inter (F^p \hat Q_i)$ contains $\hat Q_i$).  
By \lemref{Markov maps}, 
the Julia set $\hat K$ of this map is removable. All
the more, $K\subset \hat K$ is removable (one can actually see
that $K=\hat K$).
\end{pf}

\bigskip
\begin{figure}[htp]
\centerline{\psfig{figure=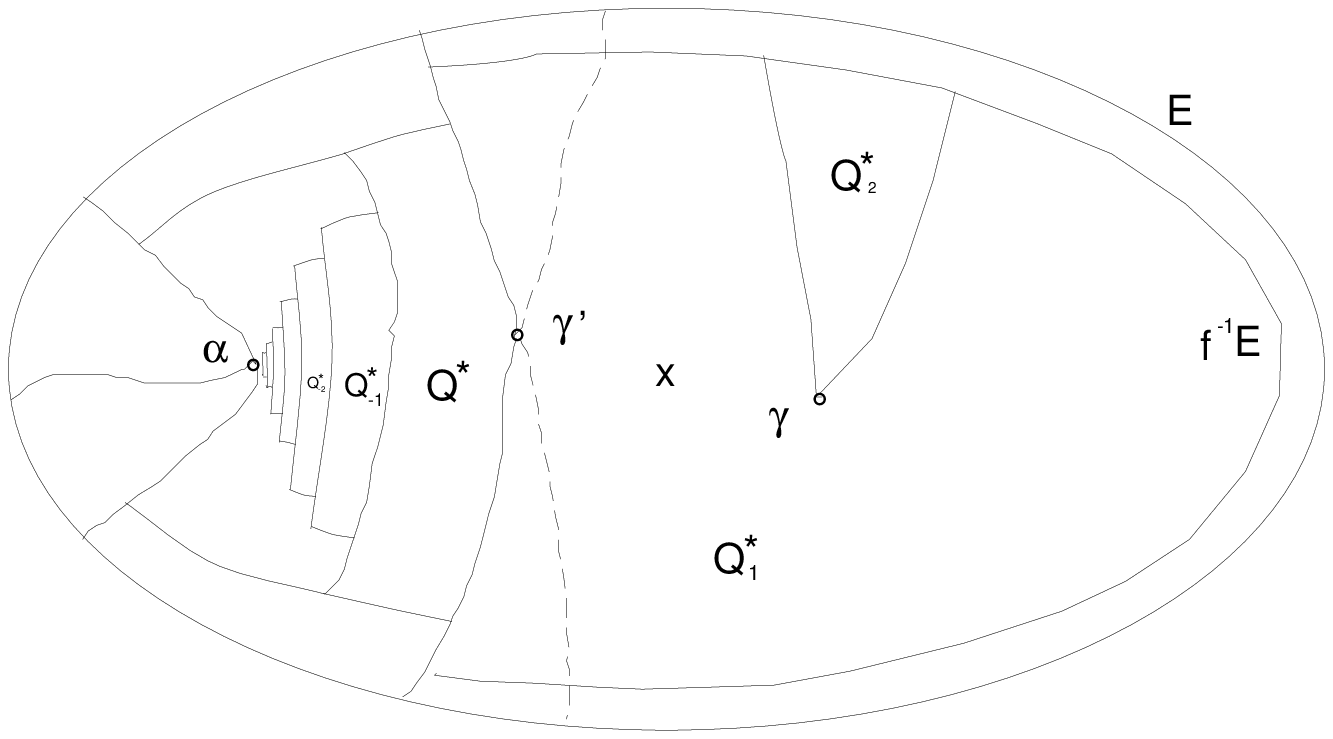,width=.9\hsize}}
\bigskip\centerline{Figure 2. Initial tiling.}\bigskip
\end{figure}

Let us now go back to \S 4.2 of Part I where
  the fundamental domain $Q$
 near the fixed point $\alpha$ is constructed. 
Recall that $\gamma\in Y^{(1)}$ is  the periodic point  of period $p$, 
$\gamma'=-\gamma$ is the
``co-periodic" point, and
$\RR(\gamma')$ is the family of rays landing at $\gamma'$.
Also, let $X=Y^{(0)}\cup_{1\leq i\leq p-1}P_i$.
 This domain is bounded by the rays landing
at $\alpha$ and equipotential $F^{-1}E$. 

 Furthermore $D$ is the connected
component of
$Y^{(1)}\sm {\RR}(\gamma')$ attached to $\alpha$,
 and $F^{-p}: D\ra F^{-p} D$ is the 
 branch of the inverse map fixing $\alpha$.

Let us also consider  quadrilaterals
$D^*=D\cap Y^{(1+p)}$ and $Q^*=Q\cap
Y^{(1+p)}$ obtained by cutting $D$ and $Q$ 
with the equipotential $F^{-p-1}E$. Note that 
$D\sm D^*=Q\sm Q^*$ consists of two quadrilaterals
 which don't contain points of the Julia set 
$J(F)$. Let $Q_{-k}^*=F^{-pk}Q^*$, $k=-1,0,1,\dots$, 
and $Q_{-2}^*=X\sm F^pD $
(see Figure 2).
Note that $J(F)\cap X$ is tiled into the pieces $Q_{-k}^*$, $k=-2,-1,\dots$

\begin{lem}\label{Q}
The hyperbolic diameter of the domains 
$Q_{-k}^*$, $k=-2,-1,0,\dots$, in $Y$ is bounded.
Moreover, if $|k-j|>1$ then there is a definite space in between $Q_{-k}$ and
 $Q_{-j}$ in $Y$.
\end{lem}  

\begin{pf} By the secondary limbs and {\it a priori} bounds assumptions,
 geometry of the configuration
 $(Y,Y^{(1)},Y^{(1+p)}, {\RR}(\gamma),  {\RR}(\gamma'))$ 
 is bounded
(see \S 4.1 of Part I).
 Hence  $Q_{-2}^*$ and $Q_{1}^*$ have a bounded hyperbolic diameter in
$Y$. For the same reason, $Q^*$ has a bounded hyperbolic diameter in $F^p D^*$.
 As $F^{-p}:  F^p  D^*\ra  F^p D^*$ is a hyperbolic contraction, 
the diameters of $ Q^*_{-k}$ in $F^p D$ are bounded by the same constant.
 All the more, they are bounded in a bigger domain $Y$. 

To prove the second statement, note that
by bounded geometry of the initial ray-equipotential configurations,
there is a definite space in between $Q_{-1}^*$ and the $Q_{-k}^*$, $k=0,1\dots$
For the same reason, 
there is a definite annulus  
$T_0\subset F^p D^*$ about $Q^*_0$ which does not intersect 
$Q_{-k}^*$, $k= 2,3\dots$. Then $T_{-i}=F^{-ip}T_0\subset F^p D^*$
 is the annulus with the same modulus
going around $Q_{-i}^*$ and disjoint from $Q_{-k}^*$ with $|k-i|>1$. 
\end{pf}

Our first essential step towards the Main Lemma is the following:

\begin{lem}\label{first step}
The Teichm\"uller distance between the configurations
 $(Y, \cup P_i, \cup Q_{-k}^*)$ 
and $(\tl Y, \cup \tl P_i, \cup\tl Q_{-k}^*)$ is bounded. 
\end{lem} 

\noindent{\bf Proof.}
 Recall that $F^s(V^0)=P_\nu$, and $F(P_i)$ univalently covers $Y$.
 Let us consider a point  $a=F^{s+1} 0\in X$. 
We will construct a qc map $(Y,a)\ra (\tl Y, \tl a)$ 
respecting the boundary marking. 

By \S 4.1 of Part I, 
geometry of the configuration $(Y,Y^{(1)},Y^{(1+p)}, {\RR}(\gamma), 
{\RR}(\gamma'))$ (and the corresponding tilde one)
is bounded , so that there is a qc pseudo-conjugacy
$$\phi: (Y,Y^{(1)}, Y^{(1+p)}, {\RR}(\gamma), {\RR}(\gamma'))
\ra (\tl Y,\tl Y^{(1)}, Y^{(1+p)}, \tl{\RR}(\gamma), \tl{\RR}(\gamma')).$$
In particular, this map conjugates $F^p:  Q^*\ra F^p Q^*$ to the corresponding
tilde map.

 As $F^p$ univalently maps $Q_{-k-1}^*$ onto $Q_{-k}^*$,
$\phi$ can be re-defined on the $Q_{-k}^*$, $k\geq 0$, 
in such a way that it becomes
the pseudo-conjugacy between the configurations
\begin{equation}\label{fist pc}
\phi: (Y,Y^{(1)},\cup Q_{-k}^*)\ra (\tl Y,\tl Y^{(1)},\cup \tl Q_{-k}^*)
\end{equation}
with the same dilatation. 
(Just let $\phi(z)=\tl F^{-kp}\circ \phi\circ F^{kp}(z)$ for $z\in Q_{-k}^*$).

It follows that $\phi(a)$ and $\tl a$ belong to the same piece of the
family $ \{Q_{-k}^*\}_{ k=-2}^\infty$ 
By \lemref{Q} the hyperbolic distance between $\phi(a)$ and $\tl a$ in 
${\tl Y}$ is bounded.

By the Moving Lemma from the Appendix,
 there is a qc map $\psi: \tl Y\ra \tl Y$
identical on the boundary and carring $\phi(a)$ to $\tl a$.
Hence $\phi_1=\psi\circ\phi: (Y, a)\ra (\tl Y, \tl a)$ is a qc map
(with a definite, though bigger, dilatation) respecting the boundary marking.

Consider now the double branched covering $F^{s+1}: (V^0,0)\ra (Y,a)$ with
the critical point at 0,
and the corresponding tilde map.
As $\phi_1: (Y,a)\ra (\tl Y, \tl a)$ 
respects the critical values for these maps,
it can be lifted to a map $\phi_2: V^0\ra \tl V^0$ with the same dilatation
respecting the boundary marking. 							
								        
Let us now construct a qc pseudo-conjugacy $\phi_3$ between corresponding
 non-critical puzzle pieces
 $P_i$ and $\tl P_i$.  It is easy as any non-central puzzle piece
$P_i$ under some iterate $F^{l_i}$ is univalently mapped onto either  $Y$
or $V^0$. In the first case let $\phi_3$ be the pullback of $\phi: Y\ra\tl Y$;
in the second let it be the pull-back of $\phi_2$. This pull-back 
preserves the dilatation and respects the boundary marking.
This provides us with a qc map $\phi_3: \cup  P_i\ra \cup\tl P_i$
respecting the boundary marking of the puzzle pieces.

 The latter property means that
$\phi_3$  matches with $h$ on $\cup\di P_i$. 
By the first part of \lemref{initial removability} and
\lemref{homotopy}, 
these maps glue together into a single homeomorphism
coinciding with $\phi_3$ on $\cup P_i$ 
and with $h$ outside, homotopic to $h$ rel $\di Y\cup\di P_i$
 (we will still denote it $\phi_3$). 

By the
Gluing Lemma from the Appendix, this homeomorphism is qc on $Y\sm K$.
By  the second part of 
\lemref{initial removability}, the residual set $K$ is removable, and
thus $\phi_3$ is automatically
 quasi-conformal  across it (with the same dilatation).
%Let us use the same notation, $\phi_3: (Y,\cup P_i)\ra (\tl Y, \cup \tl P_i)$, for the extended
%map. 
%Finally, by \lemref{initial removability} and \lemref{homotopy} below,
%\marginpar{or  Appendix?}
%$\phi_3$ is homotopic to $h$ rel $\di Y\cup\di P_i$.
\QED

The next step towards the Main Lemma is the following:

\begin{lem}\label{second step}
The configurations $(V^0,\cup V^1_i)$ and $(\tl V^0,\cup \tl V^1_i)$ 
stay bounded  
Teichm\"uller distance away.
\end{lem}

\begin{pf} Let  us consider the first return  $b=g_1 0$ 
 of the critical  point back to $V^0$. We will construct a qc map
\begin{equation}\label{first H}
H: (V^0,b)\ra (\tl V^0, \tl b)
\end{equation}
respecting the boundary marking.

  Let $u=F^{s+1} b\in X$
(where $F^s$ maps  $V^0$ onto $P_\nu$).
 Let  $\phi$ be a pseudo-conjugacy given by \lemref{first step}.Then
$\phi(u)$ and $\tl u$ belong to the same piece of the tiling 
$X\cap J(F)=\cup_{-\infty<k\leq 2} (Q_{-k}^*\cap J(F))$. 
By \lemref{Q}, the hyperbolic diameters of these pieces in $Y$ 
(and the corresponding
tilde-pieces) are bounded by a constant
$\rho$ dependent only on the selected limbs and {\it a priori} bounds.
Hence $\rho_{\tl Y}(\phi(u),\tl u)\leq \rho$.

Let $a=F^{s+1}0$, as in the proof of \lemref{first step}.
 Assume that $a\in Q_k$, $u\in Q_j$.
 Let us consider two cases:

\noindent$\bullet$ Let $|k-j|\leq 1$. Then $\rho_{\tl Y}(u,a)\leq 2\rho$. 
 Hence there is an annulus $C\subset Y$ 
going around $a$ and $u$ with $\mod C\geq \mu(\rho)>0$.
 As $F^{s+1}: (V^0,0,b)\ra (Y,a,u)$ is a
double branched covering with critical point at 0, 
the pull-back $C_0$ of this annulus to $V_0$ has
modulus at least $\mu(\rho)/2$.
 Hence $\mod (\phi(C_0))\geq K^{-1}\mu(\rho)$, where $K=\Dil(\phi)$
depends only on the selected limbs and {\it a priori} bounds. Hence
$\rho_{\tl V^0}(\phi b,0)$ is bounded.
For the same reason, $\rho_{\tl V^0}(\tl b,0)$ is bounded, and hence 
$\rho_{\tl V^0}(\phi(b),\tl b)$ is bounded.

By the Moving Lemma from the Appendix, there is a qc map 
$\psi: (\tl V^0,\phi(b))\ra (\tl V^0,\tl b)$, identical on the boundary.
 Then $H=\psi\circ\phi$ is a desired map (\ref{first H}).

\smallskip\noindent$\bullet$ Let now $|k-j|> 1$. Then by \lemref{Q},
there is a definite space in between
$Q_k^*$ and $Q_j^*$ (and between the corresponding tilde-sets).  
By the Moving Lemma, there is a qc map 
$\psi: (\tl Y,\phi(a),\phi(u))\ra (\tl Y, \tl a, \tl u)$, 
identical on $\di\tl Y$.
This map lifts to a qc map (\ref{first H}) 
(with the same dilatation).

 So, we have constructed a qc map (\ref{first H})
 which carries the critical value 
$b=g_1(0)$ to the critical value $\tl b=\tl g_1$.  \lemref{simple pull-back-2}
completes the proof. 
\end{pf}

\subsection{ Inductive  step (non-central case)}\label{sec: inductive step}

Let us now inductively estimate the Teichm\"uller  distance  between the configurations
$(V^{n-1}, \cup V^n_i)$ and $(\tl V^{n-1}, \cup \tl V^n_i)$. Let $\tau_n$
stand for the maximum of this Teichm\"uller distance and $\log \Dil(h)$,
where as above, $h$ stands for the conjugacy
 between $F$ and $\tl F$).
%\lemref{second step} provides us with the induction base ($n=1$. 
 Recall that $\mu_n=\mod(V^{n-1}\sm V^n)$ denote the principal moduli.

The following lemma is the main step of our construction.

\begin{lem}\label{inductive step} Let $\mu_n\geq\bar \mu>0$ and $\tau_n\leq \bar\tau$.
Assume that $g_n(0)\in V_k^n$ with $k\not=0$, that is, the return to level $n-1$ is non-central.
Then $\tau_{n+1} \leq \tau_n + O(\exp(-\mu_n/4))$, with a constant 
depending only on $\bar\mu$.
\end{lem}

\smallskip\noindent {\it Remark.} 
We don't assume that the non-critical  puzzle-pieces 
$V_i^n,\;i\not=0,$ are well inside $V^{n-1}$, 
since this is not the case on the levels which
immediately follow long central cascades (see Theorem II of Part I). 
%We even allow the annuli
%$V^{n-1}\sm V^n_i,\; i\not=0$ 
to be degenerate which actually happens in the beginning.

\begin{pf}  Let us skip $n$ in the notations of the objects of level $n$, so that
$V^n_i\equiv V_i$, $g_n\equiv g$, $\mu_n\equiv \mu$, etc.
Also, let $V^{n-1}\equiv\Delta$
and $g(0)\equiv c_1$. 
As above, the corresponding objects for $\tl f$ are marked with tilde.
Thus we have  two generalized polynomial-like maps $g: \cup V_i\rightarrow\Delta$ and 
$\tilde g: \cup\tilde V_i\rightarrow \tilde \Delta$, which are pseudo-conjugate by
a $K=e^\tau$-qc map 
\begin{equation}\label{zero map}
\phi: (\Delta,\cup V_i)\ra (\tl\Delta, \tl V_i).
\end{equation}
The objects on the next  level, $n+1$, will be marked with prime: $V^{n+1}\equiv V'$,
$g'\equiv g_{n+1}$  etc. (where $g'$ is {\it not} the derivative of $g$).
So $g': \cup V_j'\rightarrow \Delta'$ is the generalized renormalization
of $g$, $\Delta'\equiv V_0$. 
 
Let $\lambda(\nu)$ be the maximal hyperbolic distance between the points
in the hyperbolic plane enclosed by an annulus of modulus $\nu$.
Note that $\lambda(\nu)=O(e^{-\nu})$ as $\nu\to\infty$ 
(see Appendix A1 in Part I).
Set $\lambda=\lambda(\mu)$.
 
Our goal is to lift $\phi$ to a $K(1+O(\lambda))$-qc pseudo-conjugacy
\begin{equation}\label{final map}
\phi': (\Delta',\cup V_i')\ra (\tl\Delta', \tl V_i').
\end{equation} 
The problem is that $\phi$ need not respect the positions of the critical
values: $\phi(c)\not=\tl c$.  

\comm{
Let $\II$ be the label set for puzzle pieces $V_i$, and $\II^*=\II\sm \{0\}$.
To  any multi-index $\bar i=(i(0),i(1),\dots, i(t-1))$, $i(j)\in \II^*$ correspond
 two domains $ P_{\bar i}\subset T_{\bar i}$
such that
$$g^j P_{\bar i}\subset V_{i(j)},\; j=0,\dots,t-1,
\quad and\quad g^t P_{\bar i}=\Delta ;$$
\begin{equation}\label{T-domains}
g^j T_{\bar i}\subset V_{i(j)},\; j=0,\dots,t-1,\quad and\quad g^t T_{\bar i}=V_0
\end{equation}
(compare the proof of lemref{second step}).     end comm}

Let us consider the first landing map $T: \cup\Omega_j\ra V^0$.
By \lemref{simple pull-back},  the pseudo-conjugacy
$\phi$ admits the pull-back to a $K$-qc pseudo-conjugacy
\begin{equation}\label{next map}
\phi_1 :  (\Delta, \cup \Omega_{j})\rightarrow 
(\tilde\Delta, \cup \tilde \Omega_{j}).
\end{equation}
This {\sl localizes} the positions of the critical values  in the
sense that $\phi_1(c_1)$ and $\tl c_1$ belong to the same domain
 $\Omega_{s}\subset V_k$ (see Figure 3)
and hence the hyperbolic distance between these points
 in $\tl V_k$ is at most $\lambda$.

%\realfig{conjugacy.ps}{Localization of the critical values.}

%\bigskip
%\midinsert
\begin{figure}[htp]
\centerline{\psfig{figure=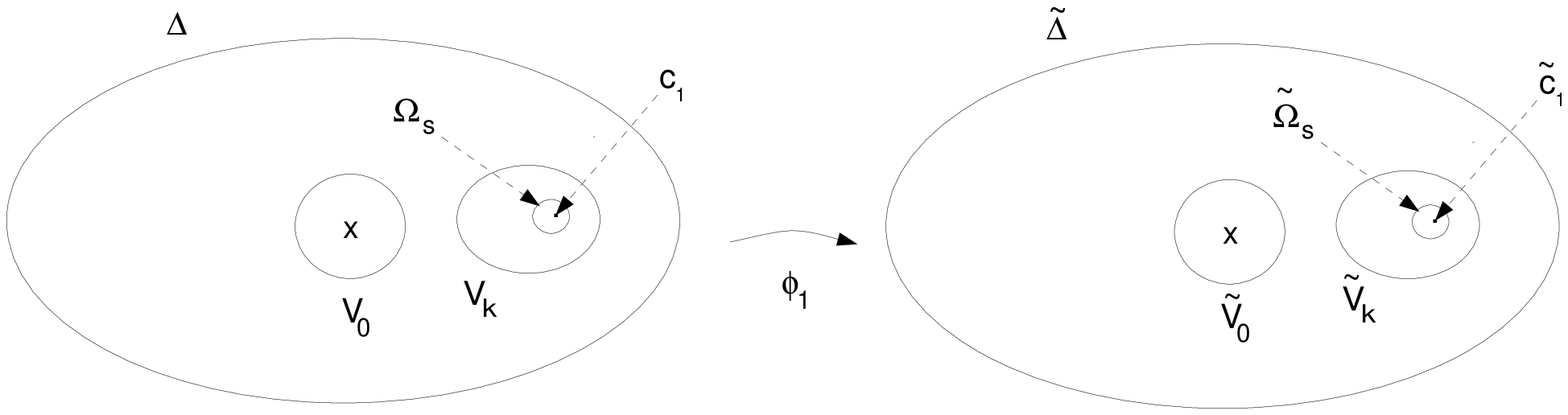,width=.95\hsize}}
\bigskip\centerline{Figure 3. Localization of the critical values.}\bigskip
%\endinsert
\end{figure}

 By the Moving Lemma from Appendix,  we can find a $(1+O(\lambda))$-qc map
\begin{equation}\label{psi}
\psi:(\tilde V_k,\phi_1(c_1))\rightarrow (\tilde V_k,\tl c_1)
\end{equation}
identical outside $\tilde V_k$.
% and moves $\phi_1 (c)$ to $\tilde c$.
Then the map 

$$\phi_2=\psi\circ \phi_1: (\Delta,\cup V_i, c)\rightarrow (\Delta, \cup V_i,\tl c)$$  is a
$K(1+O(\lambda))$-qc pseudo-conjugacy  respecting  the critical values.

Let $\{U_j'\}$ be the family of the components of the 
$\{(g|\Delta')^{-1} V_i\}$. 
The the map $\phi_2$ can be lifted to a $K(1+O(\lambda))$-qc pseudo-conjugacy
\begin{equation}\label{H}
H: (\Delta',U_i')\rightarrow (\tilde\Delta',\tilde U_i').
\end{equation}
%homotopic to $h$ rel $\Delta'\sm (\cup\inter U_i')$. 
 However $U_i'$ are not the same as $V_j'$ (components of $\{(g|\Delta')^{-1} \Omega_i\}$), 
so we have to do more:  We will localize the
positions of the critical values  $a=g'c$ and $\tilde a$ in $\Delta'$, 
and construct a   $K(1+O(\lambda))$-qc map
\begin{equation}\label{new  map}
\phi_0': (\Delta',a)\ra (\tl \Delta',\tl a)
\end{equation}
respecting the boundary marking. 
 The argument depends on the  position of $a$-points.
Let  $a_1=g(a)\in V_j$. 

\smallskip {\it Case (i).}  Assume $V_j$ is non-critical and different 
   from $V_k$. Let $a_1\in \Omega_{ l}$. Then the annuli 
  $V_j\sm \Omega_{ l}\subset V_j$ and
 $V_{k}\sm \Omega_{s}\subset V_k$ 
are disjoint (recall that $c_1\in \Omega_{s}$).
 Hence by the Moving Lemma, there is a $1+O(\lambda)$-qc map 
%\begin{equation}\label{psi1}
$$\psi_1:(\tilde \Delta,\phi_1(c_1),\phi_1(a_1))\rightarrow (\tilde\Delta,\tl c_1, \tl a_1)$$
%\end{equation}
identical outside $ (\tl V_k\cup\tl V_j)$ (where $\phi_1$ is the map (\ref{next map})). 
With  this map instead of  (\ref{psi}), the above construction leads to a map (\ref{H}) which
already  respects the critical values: $H(a)=\tl a$. Then we can let $\phi_0'=H$.
%we can simultaneously move  of $c_1$
%and $a_1$ to the right positions, and then pull the map back to $\Delta'$.

\smallskip{\sl Case (ii).} Assume that $V_j=V_k$. 

$\bullet$ Assume first that the hyperbolic diameter  of the set of four points
  $\{\tl c_1,\tl a_1, \phi_1(c_1), \phi_1(a_1)\} $ in $\tl V_k$ does not exceed
$\sqrt{\lambda}$.  
Then the hyperbolic distance between the points $\tl a_1$ and
$H(a_1)$ in
$\tl\Delta'$ is  $O(\sqrt{\lambda})$ (where $H$ is the map (\ref{H})).
 Hence there is a 
$(1+O({\lambda}^{1/4}))$-qc map
$\psi_2: (\Delta', H(a_1))\ra (\Delta',\tl a)$   identical on $\di \Delta'$. 
Define now the map (\ref{new map}) as $\psi\circ H$.

\smallskip
$\bullet$  Otherwise the hyperbolic distance between the pairs 
$(\phi_1(a_1),\tl a_1)$ 
and $(\phi_1(c_1), \tl c_1)$  in $\tl V_k$ is greater than $\sigma\sqrt{\lambda}$
(since these is an annulus of modulus $\mu$ separating one pair from another). 
Then there are separating annuli $S_i$ about these pairs with $\mod (S_i)\geq
q\sqrt{\lambda}$ (where $\sigma>0$ and $q>0$ depend only on the choice of limbs
and {\it a priori } bounds). By the
Moving Lemma, 
we can simultaneously move these points to the right positions by a
$(1+O(\sqrt{\lambda}))$-qc map 
$$\psi_2: (\tl\Delta, \tl V_k, \phi(a_1), \phi(c_1))\ra 
(\tl \Delta, \tl V_k, \tl(a_1, \tl c_1)),$$
identical on $\tl\Delta\sm \tl V_k$.
Using this map instead of (\ref{psi}) we come up with a $(1+O(\sqrt{\lambda}))$-qc map
(\ref{H}) respecting the critical values of $g$: $H(a)=\tl a$.

\smallskip{\sl Case (iii).} Let us finally assume that $V_j=V_0$ is critical.
  Then $a$ belongs to a pre-critical puzzle-piece $V'_t\subset \Delta'$.
  Since mod($\Delta'\sm V'_t)\geq \mu/2$, the hyperbolic distance between 
  $H(a)$ and $\tl a$ in $\Delta'$ is $O(\sqrt{\lambda})$ 
(where $H$ is the map (\ref{H})).
  By the Moving Lemma, there is a  $(1+O(\sqrt{\lambda}))$-qc map
  $$\psi_3: (\Delta', \phi(a))\ra  (\Delta', \tl a). $$ 
  Let us now define a map (\ref{new map}) as follows: $\phi_0'=\psi_3\circ H$.
  \smallskip 

  So in all cases we have constructed a $(1+O(\lambda^{1/4}))$-qc map (\ref{new map}).
It is still not the desired map (\ref{final map}), though.
 Now \lemref{simple pull-back-2} completes the proof. \end{pf}
\comm{ Namely, by \lemref{simple pull-back} $\phi_0'$ 
Since (\ref{new map}) matches with (\ref{zero map}) on $\di V^0$, we  can glue 
them together. Let 
$$\Phi: (\Delta, \cap V_i, a)\ra (\tl \Delta, \cup\tl V_i, \tl a)$$
coincides with $\phi$ on $\Delta\sm V^0$ and coincides with $\phi_0'$ on 
$V^0\equiv \Delta'$. By the Gluing Lemma, it is $(1+O(\lambda^{frac 1 4}))$-qc.

In the same way as (\ref{next map}) was constructed, we can turn $\Phi$ into a 
qc map (with the same dilatation)
$$\Phi_1: (\Delta, \cup P_{\bar i})\ra (\tl \Delta, \cup \tl P_{\bar i}).$$  
But now $\Phi_1(g(0))=\tl g(0)$ (since $\Phi(a)=\tl a$). Hence $\Phi_1$ admits
the lift to $\Delta'$, 
which is at last the desired map (\ref{final map}).   end comm}

\subsection{Through a central cascade} Let
$V^m\supset V^{m+1}\supset...\supset V^{m+N}$
be a cascade of central returns, so that  the critical value
$g_{m+1} 0$ belongs to $V_k,\; k=m+1,...,m+N-1$, but escapes $V^{m+N}$.

\begin{lem}\label{cascade} Let $\mu_m\geq\bar \mu>0$ and $\tau_m\leq \bar\tau$. Then
for $k\leq N+1$,
$\tau_{m+k} \leq \tau_m + O(\exp(-\mu_m/4))$,
 with a constant  depending only on $\bar\mu$.
\end{lem}

\begin{pf}
 We will adjust the proof of \lemref{inductive step} to this situation.
Let $g=g_{m+1}, \; \mu=\mod(V^m\sm V^{m+1})$, etc.
By definition, there is a $K=e^\tau$-qc pseudo-conjugacy:
$$ \phi: (V^m, \cup V^{m+1}_i)\rightarrow (\tilde V^m, \cup\tilde V^{m+1}_i).$$ 
Let us consider  the first landing map $T: \cup \Omega_j\ra V^{m+1}$
   corresponding to  $g$, $\Omega_0=V^{m+1}$. 
By \lemref{simple pull-back}, $T$ and $\tl T$ are
  pseudo-conjugate by a $K$-qc map 
$$\phi_1: (V^m, \cup \Omega_j)\ra (\tl V^m, \cup\tl\Omega_j).$$ 

Let us take a  family of puzzle pieces 
$V^{m+1}_i\subset A^{m+1}=V^m\backslash V^{m+1}$
 and pull them back to the annuli $A^{m+2},\dots,A^{m+N}.$ 
We obtain a family of puzzle pieces
$W^{m+k}_i$, together with a Bernoulli map

\begin{equation}\label{Bernoulli map}
G: V^{m+N}\bigcup_{k,i} W^{m+k}_i\ra V^m
\end{equation}
 (see \S 3.6 of Part I). Similarly let $\Omega^{m+k}_l$ stand
for the pull-backs of  the $\Omega_j\equiv\Omega^{m+1}_j$, $j\not=0$,
 to the $A^{m+k}$, $k=1,\dots N$.
If  $W^{m+k}_i\supset \Omega^{m+k}_l$ then 
$$ \mod (W^{m+k}_i\sm \Omega^{m+k}_l )\geq \mu, $$ so that
the dynamically defined  points are well localized by these puzzle pieces.

Let us now lift $\phi_1$ to the annuli
$A^{m+k}\ra \tl A^{m+k}$, $k=2,\dots, N$. We obtain a $K$-qc map
\begin{equation}
\phi_2: (V^m\sm V^{m+N},\cup  W^{m+k}_i,\cup \Omega^{m+k}_{l})\ra 
  (\tl V^m\sm \tl V^{m+N}, \cup \tl W^{m+k}_i, \cup \tl \Omega^{m+k}_{l})
\end{equation}
respecting the boundary marking.

Let $c_1\equiv g(0)\in P_l^{m+N}\subset V_k^{m+N}$.
By the Moving Lemma, there is a $(1+O(e^{-\mu}))$-qc map 
$$\psi: (\tl V^m, \tl V^{m+N}_k, \phi_2(c_1))\ra
     (\tl V^m \tl V^{m+N}_k,  \tl c_1),$$
identical outside $\tl V^{m+N}_k$. 
Then the map 
\begin{equation}\label{hole map}
\phi_3 =\psi\circ \phi_2: 
( V^m\sm V^{m+N}, \bigcup_{1\leq k\leq N,\; i\ne 0}  W^{m+k}_{i}, c_1)\ra
     (\tl V^m\sm \tl V^{m+N}, \bigcup_{1\leq k\leq N,\; i\ne 0} \tl W^{m+k}_i,  \tl c_1)
\end{equation}
is $K(1+O(e^{-mu}))$-qc, respects the
boundary marking and positions of the critical values.

Consider now the topological disks $Q_1$ and $Q_2$ in $V^{m+N}$ univalently 
mapped by $g$ onto $V^{m+N}$. The Bernoulli map $g: Q_1\cup  Q_2\ra V^{m+n}$ 
produces a family of cylinders 
$Q_{\bar i}^t$, $\bar i=(i(0),i(1),\dots, i({t-1}))$,
such that
$$g^j Q_{\bar i}^t\subset Q_{i(j)},\quad g^t Q_{\bar i}=V^{m+N}.$$
Let $\Q^t=\bigcup_{\bar i} Q^t_{\bar i}$, $Q^0\equiv V^{m+N}$.
Moreover, by \lemref{Markov maps}, the residual set 
$K=\cap\Q^t$ is removable.

The map $\phi_3$ can be consecutively lifted to the maps
$$\omega_t: \Q^{t-1}\sm \Q^t\ra \tl\Q^{t-1}\sm \tl\Q^t,\; t=1,2,\dots$$ 
with the same dilatation respecting the boundary marking.
By the Gluing Lemma, they are organized in a single qc map
 $$\omega: V^{m+N}\sm K\ra \tl V^{m+N}\sm \tl K $$ with the same dilatation.
As  $K$ is removable, this map  automatically extends across $K$:
\begin{equation}\label{H-map}
H: (V^{m+N}, \cup U_i^{m+N+1}, Q_1, Q_2)\ra 
(\tl V^{m+N},\cup \tl U_i^{m+N+1}, \tl Q_1, \tl Q_2),
\end{equation} where
$U_i^{m+N+1}\subset V^{m+N}$ are the components of $g^{-1} W^{m+N}_j$, 
$U_0^{m+N+1}\equiv V_0^{m+N+1}$.
Note that $\mod (V^{m+N}\sm U^{m+N+1})\geq \mu/2.$  
%with the same dilatation respecting the boundary marking. 

The maps (\ref{H-map}) and (\ref{hole map}) glue together into a single 
$K(1+O(e^{-mu}))$-qc map
$$\phi_4:( V^{m}, \bigcup_{1\leq k\leq N,}\bigcup_{ i\ne 0}  W^{m+k}_{i}, V^{m+N}) \ra 
        (\tl V^m, \bigcup_{1\leq k\leq N}\bigcup_{ i\ne 0} \tl V^{m+N}_i, \tl V^{m+N}).$$
Take now a family of cylinders $W^{m+k}_{\bar i}$ of the 
Bernoulli map (\ref{Bernoulli map})
 (where $\bar i$ are finite strings of symbols).
The map $\phi_4$ is naturally lifted to a qc pseudo-conjugacy $\Phi$ with the same
dilatation which respects this family of cylinders. Moreover, every
$W^{m+k}_{\bar i}$ contains a piece $V^{m+k}_{\bar i}$ such that
$$G^l V^{m+k}_{\bar i}=V^{m+k-1},\;\; where\;\; l=|\bar i|, $$
and all puzzle pieces $V^{m+k}_j$ are obtained in such a way.
  As $\phi_4$
respects the  $V^{m+k-1}$-pieces, $k\leq N$, the new map $\Phi$ respects
the $V^{m+k}_j$-pieces. Thus $\Phi$ is a  
$K(1+O(e^{-\mu}))$-qc pseudo-conjugacy  between $g_{m+k}$ and $\tl g_{m+k}$,
so that $\tau_{m+k}\leq \log K+O(e^{-\mu})$, $k=1,\dots m+N$.

\comm{*****
 As in (\ref{P-domains}),
let us consider  the family of domains
%{\sl the first landing map} 
$ \cup P^{m+k}_{\bar i}$   mapped  onto $V^{m+N}$ by appropriate iterates  of $G$.

Let us proceed further along the lines of the proof of \lemref{inductive step}
using  the following substitution: the $\{W^{m+k}_i\}$ play the role of
the non-central pieces $\{V_i\}$, $V^{m+N}$ plays the role of $\Delta'$.
So we pull $h$ back to 
$$h_1: (\Delta, \cup P_l, V^{m+N})\rightarrow
  (\tilde\Delta, \tilde\cup P_l, \tilde V^{m+N}), $$
correct this map to make it  respect the $g$-critical values and then
pull it back to $V^{m+N}$ as in (25):
$$H: (V^{m+N}, U_i^{m+N+1})\rightarrow (\tilde V^{m+N}, \tilde U_i^{m+N+1}), $$
where $U_i^{m+N+1}$ are the pull-backs of $W^{m+N}_i$ and $V^{m+N}$.

Let us proceed further along the lines of the proof of \lemref{inductive step}
 using  the $\{W^{m+N}_i\}$  in place of
$\{V_i\}$ and  $V^{m+N}$  in place of $\Delta'$. *********}

Let us proceed with the estimate of $\tau_{m+N+1}$.
Take the first return $a$ of the
critical point back to $V^{m+N}$,  and construct a $K(1+O(e^{-\frac \mu 4}))$-qc  map
\begin{equation}\label{map respecting a}
\phi_0': (V^{m+N},a)\ra (\tl V^{m+N}, \tl a)
\end{equation}
To this end let us
go  through Cases (i), (ii),
(iii) of the proof of \lemref{inductive step} using  the $\{W^{m+N}_i\}$  in place of
$\{V_i\}$ and  $V^{m+N}$  in place of $V^{m+1}\equiv\Delta'$.

In the first two cases the argument is the same  as above. 
However, the last case is
different since the pre-critical puzzle-pieces $Q_1$ and $Q_2$
are not necessarily well inside of $V^{m+N}$. 
To take care of this problem
let us  consider the first ''escaping moment" $t$ when 
$b\equiv g^t a\not\in Q_1\cup Q_2$. Then $b\in U^{m+N+1}_i$ for some $U$-domain
from (\ref{H-map}). 
Then there is a domain $\Lambda\subset Q_1\cup Q_2$ containing $a$
which is univalently mapped onto $U^{m+N+1}_i$ by $g^t$. Moreover
$$\mod(Q\sm \Lambda)=\mod(V^{m+N}\sm U^{m+N+1}_i)\geq \mu.$$
By means of $g: Q_1\cup Q_2\ra V^{m+N}$,
the map (\ref{H-map}) can be turned into a  qc map (with the same dilatation) 
$$H_1: (V^{m+N}, \Lambda)\ra (\tl V^{m+N}, \tl \Lambda)$$
(coinciding with $H$ outside $Q_1\cup Q_2$). 
This gives us an appropriate localization of the $a$-points.
The Moving Lemma  turns $H_1$ into (\ref{map respecting a}).

\lemref{simple pull-back-2} completes the proof.   \end{pf}

\comm{Finally, let us  turn ({map respecting a}) into the desired pseudo-conjugacy
$$\phi': (V^{m+N}, \cup V^{m+N+1}_i)\ra (\tl V^{m+N}, \cup \tl V^{m+N+1}_i)$$
with the same dilatation. 
The argument is identical with the final part of the proof of 
\lemref{inductive step}, with usage of the Bernoulli map (\ref{Bernoulli map})
in place of $g$. This gives a desired estimate for $\tau_{m+N+1}$.  end comm}

\subsection{Proof of the Main Lemma}
Let $\{i(k)\}$ be the sequence of non-central levels in the principal nest
$V^0\supset V^1\supset\dots$  Let  $i(n-1)+1<m\leq i(n)+1.$
By  \lemref{cascade}, 
%the Teichm\"uller distance $\tau_n$ between $(V^{n-1},
%V^n)$  and $(\tl V^{n-1}, \tl V^n)$ is bounded by
\begin{equation}\label{sum}
\tau_{m}\leq\log K^*+O(\sum_{k=0}^{n-1}\exp(-{\frac 1 4} \mu_{i(k)+1})).
\end{equation}
 But by Theorem III from Part I  \cite{L6},  
the principal moduli
$\mu_{i(k)+1}$ grow at linear rate: 
$\mu_{i(k)+1}\geq Bk,$ where the constant $B$ depends only on $\mu_1$. 
Hence  the sum (\ref{sum}) is bounded by $\log K_*+C(\mu_1)$. 

In turn, by Theorem I of Part I
the modulus $\mu_1$ is bounded by a constant depending only 
on the selected limbs and
{\it a priori} bounds. Hence 
$\tau_n\leq \log K_* + A$, where $A$ depends only on the choice of limbs and 
{\it a priori} bounds. 
 The Main Lemma is proved.  \QED

\subsection{Last cascade}\label{sec:last cascade}
 If the map $F\equiv F_m=R^m f$ is not renormalizable
then  the principal nest consists of infinitely many central cascades, and
 the Main Lemma gives a uniform bound on the Teichm\"uller distance between
the corresponding generalized renormalizations.

Otherwise the principal nest ends up with an infinite central cascade
$V^{n-1}\supset V^n\supset\dots$ shrinking to the little Julia set $J^{m+1}$
of the next renormalization $g_n = F_{m+1}\equiv R^{m+1} f$.
 All levels $ n, n+1, \dots$  of this final cascade
 are called the  {\it renormalization levels}. 

\begin{lem}\label{last cascade}
Let $n$ be  a renormalization level and 
$H: (V^{n-1}, V^n)\ra (\tl V^{n-1}, \tl V^n)$ be a $K$-qc pseudo-conjugacy
between $g_n$ and $\tl g_n$. Then there is a homeomorphism 
$\phi: (V^{n-1}, J^{m+1})\ra (\tl V^{n-1}, \tl J^{m+1})$ homotopic to $h$
rel $(J^{m+1}\cup \di V^{n-1})$, and $K$-qc on $V^{n-1}\sm J^{m+1}$.
\end{lem}

\begin{pf}
 Recall that $k^n=V^{k-1}\sm V^k$.
The map  $H: A^{n}\ra \tl A^{n}$  admits a lift
to qc maps (with the same dilatation) 
$H_k: A^{n+k}\ra \tl A^{n+k}$   homotopic to $h$ rel the annuli boundary.
These maps match to a single qc map 
$\phi: V^{n-1}\sm J^{m+1}\ra \tl V^{n-1}\sm \tl J^{m+1}$
with the same dilatation conjugating $F_{m+1}$ to $\tl F_{m+1}$. 
By \corref{matching}, this map  (and the whole homotopy between it and $h$)
matches with $h$ on $J(F_{m+1})$.
\end{pf}

\subsection{Spreading around}\label{spreading}
Let us consider the pieces $P_j\subset Y\equiv Y^{(0)}$ of the initial partition 
(\ref{initial partition}), and the Markov map
$G: \cup P_i\ra Y$ (see \S 3.2 of Part I). 
Let us consider the {\it first landing map} to $V^0$,
 $T_0: \cup \Omega_i^0\ra P_0.$
By \lemref{first step} and \lemref{simple pull-back},
there is a qc pseudo-conjugacy
$\phi_0: (Y, \cup \Omega^0_i)\ra (\tl Y, \cup \tl \Omega^0_i).$
Let us also consider the following maps:

\smallskip\noindent$\bullet$ The first landing  maps to $V^n$ corresponding to the
generalized renormalization $g_n: \cup V^n_i\ra V^{n-1}$:
$$T_n: \cup \Omega^n_i\ra V^{n},\;\quad  \Omega^n_i\subset V^{n-1};$$

\smallskip\noindent$\bullet$ 
The first landing  maps to $V^n$ corresponding to $G$:
$$S_n: \cup O^n_i\ra V^{n},\;\quad  O^n_i\subset Y.$$
Clearly  
 \begin{equation}\label{S-T} 
S_0=T_0\;\; and \;\;   S_n=T_n\circ S_{n-1}.
\end{equation}
By the Main Lemma and \lemref{simple pull-back}, there is a sequence of qc 
pseudo-conjugacies 
$$\phi_n: (V^{n-1}, \cup \Omega^n_i)\ra (\tl V^{n-1}, \cup \tl \Omega^n_i),
\quad n<N+1,$$
where $N$ is the first  DH renormalizable level
 (if $F$ is non-renormalizable then
$N=\infty$). Let us turn it inductively into a sequence of pseudo-conjugacies
\begin{equation}\label{H-sub-n}
 H_n: (Y, \cup O^n_i)\ra (\tl Y, \cup \tl O^n_i)
\end{equation}
between $S_n$ and $\tl S_n$ (with the same dilatation).
Indeed, using (\ref{S-T}), we can define it as  follows:
$$H_n| O_i^{n-1}=  (\tl S_{n-1}|\tl O^{n-1}_i)^{-1}\circ (\phi_n|V^{n-1})\circ
S_{n-1}|O^{n-1}_i.$$
As these maps match with $H_{n-1}$ on the boundaries $\di O^{n-1}_i$, the glue together
into  single qc conjugacy (\ref{H-sub-n}) with the same dilatation.

If $F$ is non-renormalizable, we obtain an  infinite sequence of qc pseudo-conjugacies
$H_n$ (with uniformly bounded dilatation).
As the pieces $V^n_i$ shrink as $n\to\infty$, there is the limit qc  map 
\begin{equation}\label{limit map}
H: (Y, J(F)\cap Y)  \ra (\tl Y, \tl J(F)\cap \tl Y) 
\end{equation}
homotopic to $h: J(F)\cap Y\ra \tl J(F)\cap \tl Y$ rel $\di Y\cup J(F)$.

\comm{Let $\Theta$ be the component of the standard neighborhood $U\supset J(F)$...
  
Doing the same construction on every puzzle piece $ Y^{(0)}_i$, $i=1,\dots p-1$
(where $p$ is the number of rays landing at $\alpha$-fixed point of $F$),
we obtain a qc map } 

Assume $F$ is renormalizable. 
Let $\II$ be the family of little Julia sets $J^{m+1}_i$
contained in $Y$, $J^{m+1}\equiv J^{m+1}_0$. Let us consider the 
last pseudo-conjugacy  (\ref{H-sub-n}) on the renormalization level $N$.
Let us replace it on $V^{n-1} $ by the pseudo-conjugacy 
$$\phi_N: (V^N, J^{m+1})\ra (\tl V^N, \tl J^{m+1})$$
 constructed in \lemref{last cascade}.
Spread it around by the landing map $S_{N}$,  that is, set
$$H|O_N=     (\tl S_N|\tl O_N)^{-1}  (\phi_N|V_N)\circ S_N|O_N.$$
As these maps match on the  $\di O_N$ with $H_N$, 
they glue together into a homeomorphism
\begin{equation}\label{limit map-1}
H: (Y, \bigcup_{i\in \II} J^{m+1}_i)\ra
 (\tl Y, \bigcup_{i\in\II} \tl J^{m+1}_i),
\end{equation}
quasi-conformal on $Y\sm \bigcup_{i\in \II}J^{m+1}_i$ 
(with dilatation depending only on
the  choice of limbs and {\it a priori} bounds), and homotopic to $h$
rel $\di Y\bigcup_{i\in \II} J^{m+1}_i$. 

%Given a puzzle piece $V^n_i$, $n\geq 1$,  let $t(n,i)$ stand for
%its return time, that is, $F^{t(n,i)} V^n_i=V^{n-1}$. Let $\VV$ be the
%whole family of pieces $F^t V^n_i$, $t=1,\dots , t(n,i)$, for all $V^n_i$.
%Clearly  any two pieces of this
%family are either disjoint or coincide. Let $\V=\bigcup_{X\in \VV} X$.

%Let us assume first that the Julia sets $J^m_j$ don't touch,
%that is, $F_{m-1}$ is not immediately renormalizable. Then
 Let us
consider the backward orbit  $Y\equiv Y_0, Y_{-1},\dots, Y_{-r+1}$ 
of $Y$ under $f$ 
such that $Y_{-k}\ni f^{r-k}0$, where $r$ is the first return time of
the critical orbit to $Y$. The disks $Y_{-k}$ have disjoint interiors.
Let us pull the map $H$ back to these
disks, that is, set 
$$h_{m+1}| Y_{-k}=  (\tl f^k|\tl Y_{-k})^{-1}  H\circ   f^k|Y_{-k}.$$
%This map conjugates $f|\J^{m+1}$ to $\tl f|\tl \J^{m+1}$,  and its
%dilatation  depends only on the choice of limbs and
%{\it a priori} bounds.
As this map respects the boundary marking of the $Y_{-k}$, it extends to
to the whole plane as $h_m$, which provides the desired next approximation
to the Thurston conjugacy (see \secref{strategy}). 

The Rigidity Theorem is proved.

\comm{
Assume now that $F_{m-1}$ is immediately renormalizable.
Let us consider the external rays 
of $F_{m-1}$ landing at the fixed point $\alpha_{m-1}$
Let $Y_\#$ be the part of $Y$ bounded by the external rays of $F_{m-1}$
landing at the fixed points of $J^m$ 
Note that the limit map $H$ from (\ref{limit map}) or (\ref{limit map-1})
respects the rays landing at the fixed points of $J(F)$. Indeed, the initial 
map (\ref{initial map}) does, and all further
 pull-backs preserve this property.
 
Doing the same construction  on every 
puzzle piece $ Y^{(0)}_i$, $i=1,\dots p-1,$
(where $p$ is the number of rays landing at $\alpha$-fixed point of $F_m$),
we obtain a  pseudo conjugacy  near the  Julia set $J^m$. 
Doing this for  every little Julia set $J^{m}_j$
we obtain a desired map 
$$h_{m+1}: (\U^m,\J^{m+1})\ra (\tl \U^m, \tl \J^{m+1})$$
(see \secref{strategy}). 
This completes the proof of the Rigidity Theorem.}

\section{Appendix: Quasi-conformal maps}\label{appendix}

\subsection{}  There are 
a few Russian and English sources  on the basic theory of quasi-conformal
 maps:  \cite{A,B,Kr,LV,V}.

 A homeomorphism $h: U\rightarrow V$,
where $U, V\subset \C$,
 is called quasi-conformal (qc) if it has locally integrable
distributional derivatives $\di h$, $\bar\di h$, and  
$|\bar\di h/ \di h|\leq k<1$ almost everywhere. 
As  this local
definition is conformally invariant, one  can define  qc 
homeomorphisms between Riemann surfaces.

 One  can associate to a qc map an analytic object called Beltrami
differential, namely 
$$\mu={\bar\di h\over \di h} {d\bar z\over dz},$$ with
$\|\mu\|_{\infty}<1$. The corresponding  geometric object is  a
measurable  family of infinitesimal ellipses (defined up to dilation), 
pull-backs by $h_*$ of the field of infinitesimal circles. The
eccentricities of these ellipses are ruled by $|\mu|$,
 and are uniformly bounded almost everywhere, while the orientation of
the ellipses is ruled by the $\arg\mu$. The dilatation
$\Dil(h)\equiv K_h=(1+\|\mu\|_{\infty})/(1-\|\mu\|_{\infty})$  of $h$ is the
essential supremum of the eccentricities of these ellipses. 
 A qc map $h$
is called $K$-qc if $\Dil(h)\leq K$. 

\proclaim Weil's Lemma. A 1-qc map is analytic.

One of the most remarkable facts of analysis is that 
{\it any Beltrami differential with
$\|\mu_{\infty}\|<1$  (or rather a measurable field of ellipses with
essentially bounded eccentricities) is locally generated by a qc map},
unique up to post-composition with an analytic map. 
 Thus such a Beltrami differential on a Riemann surface $S$  induces a
conformal structure quasi-conformally equivalent to the original
structure of $S$,    Together with the Riemann mapping theorem this
leads to the following statement:

\proclaim Measurable Riemann Mapping Theorem.
 Let $\mu$ be a  Beltrami differential on
$\bar\C$ with  $\|\mu_{\infty}\|<1$, Then there is a quasi-conformal
map
$h: \bar\C\rightarrow \bar\C$ which solves the Beltrami equation:
$|\bar\di h/ \di h|=\mu$.

In what follows  by a conformal structure we will mean a structure
associated to measurable Beltrami differentials $\mu$ with
$\|\mu\|_\infty<1$.  We will denote  by $\sigma$ the standard
structure corresponding to zero  Beltrami differential.

Another fundamental property of the space of qc maps is compactness:

\proclaim Compactness Lemma. The space of $K$-qc maps $h:
\C\rightarrow \C$ normalized by $h(0)=0$ and
$h(1)=1$ is compact in the uniform topology on the Riemann sphere.

The following gluing property is also important:

\proclaim Gluing Lemma. Let us have two disjoint domains $D_1$ and
$D_2$ with a piecewise smooth  arc $\gamma$ of their common boundary. Let
$D=D_1\cup D_2\cup \gamma$. If $h: D\rightarrow \C$ is a homeomorphism
such that $h|D_i$ is 
$K$-qc, then $h$ is $K$-qc.

One of Sullivan's leading ideas  was the idea of the Teichm\"uller
metric on the space of deformations of a conformal dynamical systems.
The prototype for this metric is the classical Teichm\"uller metric on
the space of marked Riemann surfaces. 
 A  Riemann surface (perhaps with boundary)
is said to be
{\it marked} if it is endowed with 
 a preferred basis of the fundamental group and
a preferred parametrization  of the boundary components. 
The Teichm\"uller distance
$\dist(S_1, S_2)$ between two marked Riemann surfaces  
is defined as the infimum of the
dilatations
$K_h$, where $h : S_1\rightarrow S_2$ runs over qc homeomorphisms in
the homotopy class respecting the marking.  

Let $D$ be a simply connected domain conformally 
equivalent to the hyperbolic plane $\bH^2$.
Given a family of subsets $\{S_k\}_{k=1}^n$ in $D$,
 let us say that a family of disjoint
annuli
$A_k\subset D\sm \cup S_k$ is {\it separating} if $A_k$ surrounds $S_k$ but does not
surround the $S_i$, $i\not= k$.
The following lemma is used in the present paper uncountably many times:

\medskip\noindent{\bf Moving Lemma.}  {\it  
$\bullet$ Let  $a, b\in D$ be two points on hyperbolic distance  $\rho\leq \bar\rho$. Then
there is a diffeomorphism
$\phi:(D,a)\ra (D,b)$, 
identical near $\di D$, with dilatation $\Dil(\phi)=1+O(\rho)$,
where the constant depends only on $\bar\rho$.

$\bullet$ Let $\{(a_k,b_k)\}_{k=1}^n$ 
be a family of  pairs of points which admits
 a family of separating annuli $A_k$ with  $\mod A_k\geq\mu$. 
Then there is a diffeomorphism
$\phi:(D,a_1,\dots a_n)\ra (D,b_1,\dots, b_n)$,
 identical near $\di D$, with dilatation
$\Dil(\phi)=1+O(e^{-\mu})$.}

\begin{pf} As the statement is conformally equivalent, we can work  with the 
unit disk model of the hyperbolic plane, and  can also assume that $a=0$. 
Also, it is enough to prove the statement for sufficiently small $\rho$.

There is a smooth function $\psi: [0,1]\ra [\rho,1]$ 
such that  $\psi(x)\equiv \rho$
near 0, $\psi(x)\equiv 0$ near 1, and $\psi'(x)=O(\rho)$, 
with  a constant depending
only on $\bar \rho$.

Let us define a smooth map $\phi: (\D,0)\ra (\D,b)$ as $z\mapsto z+\psi(|z|)$.
Then 
\begin{equation}\label{derivatives}
\di\phi(z)=1+\psi'(|z|){\bar z\over 2|z|}=1+O(\rho),
\quad  \bar\di\phi(z)=\psi'(|z|){z\over 2|z|}=O(\rho).
\end{equation}
Thus
$$ Jac(f)=|di\phi(z)|^2-|\bar\di\phi(z)|^2=1-O(\rho).$$
Hence for sufficiently small $\rho>0$,
 $f$ is a local orientation preserving diffeomorphism.
As $f: \di \D\ra \di \D$, $f$ is a proper map. Hence it is a diffeomorphism.

Finally, (\ref{derivatives}) yields that the Beltrami coefficient 
$\mu_f=O(\rho)$, so that
 the dilatation $\Dil(f)=1+O(\rho)$.
\end{pf}

Let $Q\subset\C$, $h: Q\ra \C$ be a homeomorphism onto its image.
It is called {\it quasi-symmetric  (qs)} if for any three points
$a,b,c\in Q$ such that
$$ q^{-1}\leq {|a-b|\over |b-c|}\leq q,$$
we have:
$$ \kappa(q)^{-1}\leq {|a-b|\over |b-c|}\leq \kappa(q).$$
It is called $\kappa$-quasi-symmetric if $\kappa(1)\leq\kappa$.
It follows from the Compactness Lemma that any $K$-qc map is 
$\kappa$-quasi-symmetric, with a $\kappa$ depending only on $K$. 

\comm{
Let us discuss quasi-symmetric maps of the circle $\T=\{z: |z|=1\}$. 
 Given an interval
$J\subset \T$, let $|J|$ denote its length.  An orientation preserving
 map $h: \T\ra \T$ is called
$\kappa$-{\it quasi-symmetric} ($\kappa$-qs) 
if for any two adjacent intervals $I,J\subset\T$,
$|h I|/|h J|\leq \kappa$. 

Let $\T_r=\{z: |z|=r\}$, $\T\equiv \T_1$. 
Let $\A(r, R)=\{z: r<|z|<R\}$. Similar notations are used for a closed annulus
$\A[r,R]$ (or semi-closed one).}

\proclaim Ahlfors-B\"orling Extension Theorem. Any 
$\kappa$-quasi-symmetric   
map $h: \T\ra T$   
 extends to a $K(\kappa)$-qc map $H: \C\ra \C$. 
Vice versa: The restriction of any
%normalized 
%$K$-qc map $H: (\A(r^{-1}, r),\T,1,-1,i)\ra (U, \T, 1,-1,i))$,
% where $U\subset \C$, to the
$K$-qc map 
$H: (\A(r^{-1}, r),\T)\ra (U, \T)$ (where $U\subset \C$) to the
circle  $\kappa(K,r)$-quasi-symmetric. 

Let us note that in the upper half-plane model, the Ahlfors-B\"orling extension
of a qs map $\R\ra \R$ is affinely equivariant (that is, commutes with
the action of the complex affine group $z\mapsto az+b$).

%Note that the second part is false without normalization. Indeed, M\"obius maps
%$M: (\C,\T)\ra (\C,\T)$ being conformal don't have any bound on the quasi-symmetric
%dilatation on the circle. However, this is the only problem: Any $K$-qc map
%$H: (\C,\T)\ra (\C,\T)$ is equal to $M\circ h$ where $M$ is a M\'obius map
%preserving the circle and  $h|\T$ is $\kappa(K)$-quasi-symmetric. 

%For us a different normalization will be useful. 

%\proclaim Corollary. Let $H: (A(r^{-1}, r),\T)\ra (U, \T)$ be a $K$-qc map commuting with
%the central reflection $z\mapsto -z$.  Then $H|\T$ is $\kappa(K,r)$-quasi-symmetric.

%\begin{pf} We can normalize this map by a rotation so that $H(-1)=-1,\; H(1)=(1)$.
%Let us show that $\dist(H(i), \{1,-1\})\geq \delta(K,r)$. Indeed, there is 
%an annulus $L\subset \A(r^{-1}, r)\sm \{1,-1, i,-i\}$ separating $\{1,-i\}$ from
%$\{-1, i\}$, with modulus $\geq \mu(r)$. Hence $\mod (H(L))\geq K^{-1}\delta(r)$.
%This implies that 
%\end{pf}

\proclaim Interpolation Lemma. Let  us consider two round annuli $A=\A[1,r]$ and
$\tl A=\A[1,\tl r]$, with $0<\eps\leq \mod A\leq \eps^{-1}$ and
 $\eps\leq \mod \tl A\leq\eps^{-1}$.
 Then  any $\kappa$-qs map 
$h: (\T, \T_r)\ra  (\tl \T, \tl \T_{\tl r})$ admits a 
$K(\kappa, \eps)$-qc extension to a map $H: A\ra \tl A$. 

\begin{pf} Since $A$ and $\tl A$ are $\eps^2$-qc equivalent, we can assume
without loss of generality that $A=\tl A$. Let us cover $A$ by the upper half-plane,
$\theta: \bH\ra A$, $\theta(z)=z^{-\log r i\over \pi}$, where the covering group generated
by the dilation $T: z\mapsto \lambda z$, with $\lambda=e^{2\pi^2\over \log r}$. 
Let $\bar h: (\R,0)\ra (\R,0)$ be the lift of $h$ to $\R$ such that 
$\bar h(1)\in [1, \lambda)\equiv I_\lambda$ and  $\bar h(1)\in ( -\lambda, -1]$ 
(note that
$\R_+$ covers $\T_r$, while $\R_-$ covers  $\T$). Moreover, since $\deg h=1$,
it commutes with the deck transformation $T$.

A  direct calculation shows that the dilatation of the covering map $\theta$
 on the fundamental intervals $I_\lambda$ and $-I_\lambda$ is comparable with
$(\log r)^{-1}$. Hence $\bar h$ is $C(\kappa, r)$-qs on this interval. By equivariance it
is $C(\kappa, r)$-qc on the rays $\R_+$ and $\R_-$. 

It is also quasi-symmetric near the origin. Indeed, by the equivariance and normalization,
$$(1+\lambda)^{-1}|J|\leq |\bar h(J)|\leq (1+\lambda) |J|$$
for any interval $J$ containing 0, which easily implies quasi-symmetry.

Since the Ahlfors-B\"orling extension is affinely equivariant,
 the map $\bar h$ extends to a $K(\kappa, r)$-qc map
$\bar H: \bH\ra \bH$ commuting with $T$. Hence $\bar H$ descends to a  
$K(\kappa, r)$-qc map $H: A\ra A$.
\end{pf}

\subsection{Removability}\label{sec: removability}

A compact set $X\subset \C$ is called {\it removable} if for any neighborhood
$U\supset X$, any conformal map $h: U\sm X\ra \C$ 
admits a conformal extension across $X$.  
Let us show that removability is quasi-conformally invariant.

%Given an open set $V\subset \C$, let us call a point $z\in \C\sm X$ {\it essential}
%if there exists a Jordan curve $\gamma\subset V$ enclosing $z$. Let us call the set of
%essential points the {\it essential complement} of $V$.

%A compact set $X\subset \C$ is said to have an {\it absolute measure zero} if 
%for any conformal map $h: U\sm X\ra \C$ as above, the essential complement of
%$h (U\sm X)$ has zero Lebesgue measure.

%\proclaim Removability Criterion [see {SN}]. A compact set $X$ is removable if and 
%only if it has absolute measure zero.

\begin{lem}\label{qc invariance}
Let $\phi: (\C,X)\ra (\C, \tl X)$ be a qc map.
If the set  $X$ is removable then $\tl X$ is removable as well.
\end{lem} 

\begin{pf} Let $\sigma$ be the standard conformal structure on $\C$.
Let $\tl U\supset \tl X$ be a neighborhood of $\tl X$, and let 
$\tl h: \tl U\sm \tl X\ra \C$ be a conformal map. Let as consider a conformal 
structure $\tl \mu$ on $\C$ which is equal to 
$(h\circ\phi)_*(\sigma)$ on $h(\tl U\sm \tl X)$, and
is equal to $\sigma$ outside.  
By the Measurable Riemann Mapping Theorem,
there is a qc map $\psi: \C\ra \C$ such that $\tl \mu=\psi_*(\sigma)$.

Let $U=\phi^{-1}\tl U$. Then the function 
$h=\psi^{-1}\circ \tl h\circ \phi: U\sm X\ra \C$
is conformal. As $X$ is removable, it admits a conformal extension across $X$.
We will use the same notation $h$ for the extended function. Then the formula
$\tl h=\psi\circ h\circ\phi^{-1}$
 provides us with a conformal extension of $\tl h$
across $\tl X$.
\end{pf}

Let us now show that removable sets are also qc-removable.

\begin{lem}\label{qc removability}
 Let $X$ be a removable set and $U\supset X$ be
its neighborhood. Then any qc map $h$ on $U\sm X$ admits a qc 
extension across $X$.
\end{lem}

\begin{pf} Let us consider a conformal structure $\mu$
 equal to $h^*(\sigma)$ on $U\sm X$,
and equal to $\sigma$ on the rest of $\C$. 
By the Measurable Riemann Mapping Theorem,
there exists a qc map $ \phi:\C\ra\C$ such that $\mu=\phi^*(\sigma)$. 
Then the function $\tl h= h\circ \phi^{-1}$ is univalent on
 $\tl U\sm \tl X\equiv \phi U\sm\phi X$.

 By \lemref{qc invariance}, the set $\tl X$
is removable.  Hence $\tl h$ admits a conformal extension across $\tl X$.
Then the formula 
$h=\tl h\circ \phi$ provides us with a qc extension of $h$ across $X$.
\end{pf}

Let us finally state a simple 
condition for removability (see, e.g., [SN])
 which is used many times in this paper.

\proclaim Removability Condition.
 Let $X$ be a Cantor set satisfying the following
property. There is an $\eta>o$ such that for any point $z\in X$, 
there is a nest
of disjoint annuli $A_i(z)\subset \C\sm X$ surrounding $z$ with
$\mod(A_i(z))\geq \eta$.  Then $X$ is removable.

\comm{*********************
    Index
Part I.  Intro.
Boundary marking \S 2.6
Secondary limb condition.  p.1.
$\A(r,R)$ annulus
$\T_r$ circle, $\T$
$\D_r$ disk, $\D$
*************************}

\end{document}